\documentclass[11pt]{amsart}
\usepackage[latin1]{inputenc}
\usepackage{amsfonts}
\usepackage{amsthm}
\usepackage{amsmath}
\usepackage{amssymb}
\usepackage{amscd}
\usepackage{verbatim}
\usepackage{multicol}
\usepackage{tabularx}

\theoremstyle{definition}
\newtheorem{thm}{Theorem}[section]
\newtheorem{prop}[thm]{Proposition}
\newtheorem{cor}[thm]{Corollary}
\newtheorem{con}[thm]{Conjecture}
\newtheorem{lem}[thm]{Lemma}
\newtheorem{defn}[thm]{Definition}

\newtheorem*{rem}{Remark}

\numberwithin{equation}{section}

\def\cal#1{\text{$\mathcal{#1}$}}

\def\ord#1^#2{#1$^{\text{#2}}$}

\def\lie#1{\mathfrak{#1}}
\def\tlie#1{\tilde{\mathfrak{#1}}}
\def\hlie#1{\hat{\mathfrak{#1}}}

\def\uqr#1^#2{\text{$U_q^{#2}(\lie #1)$}}

\def\uqhr#1^#2{\text{$U_q^{#2}(\hlie #1)$}}
\def\us#1^#2{\text{$U_{\xi}^{#2}(\lie #1)$}}
\def\ush#1^#2{\text{$U_{\xi}^{#2}(\hlie #1)$}}
\def\dus#1^#2{\text{$\dot{U}_{\xi}^{#2}(\lie #1)$}}
\def\dush#1^#2{\text{$\dot{U}_{\xi}^{#2}(\hlie #1)$}}
\def\gb#1{{\mbox{\boldmath $#1$}}}

\def\wt{{\rm wt}}

\def\opl_#1^#2{\text{\scriptsize$\bigoplus\limits_{\text{\footnotesize$#1$}}^{\text{\footnotesize$#2$}}$}}
\def\oti_#1^#2{\text{\scriptsize$\bigotimes\limits_{\text{\footnotesize$#1$}}^{\text{\footnotesize$#2$}}$}}

\renewcommand{\thefootnote}

\begin{document}

\flushbottom

\title[Finite-Dimensional Representations of Hyper Loop Algebras]{Finite-Dimensional Representations of\\ Hyper Loop Algebras}

\author[D. Jakeli\'c and A. Moura]{Dijana Jakeli\'c and Adriano Moura}
\address{Max-Planck-Institut f\"ur Mathematik, D-53111, Bonn, Germany} \email{jakelic@mpim-bonn.mpg.de}

\address{UNICAMP - IMECC, Campinas - SP, 13083-970, Brazil.} \email{aamoura@ime.unicamp.br}

\maketitle
\centerline{\small{
\begin{minipage}{350pt}
{\bf Abstract:} We study finite-dimensional representations of hyper loop algebras, i.e.,  the
hyperalgebras over an algebraically closed field of positive characteristic associated to the loop algebra over a complex finite-dimensional simple  Lie algebra. The main results are the
classification of the irreducible modules, a version of Steinberg's
Tensor Product Theorem, and the construction of positive
characteristic analogues of the Weyl modules as defined by Chari and
Pressley in the characteristic zero setting.  Furthermore, we start the study of reduction modulo $p$ and prove
that every irreducible module of a hyper loop algebra can be constructed as a quotient of a
module obtained by a certain reduction modulo $p$ process applied to a suitable characteristic zero module. We conjecture that the Weyl modules are also obtained by reduction modulo $p$. The conjecture implies a tensor product decomposition for the Weyl modules which we use to describe the blocks of the underlying abelian category.
\end{minipage}}}

\section*{Introduction}

Let $G$ be a semisimple connected algebraic group over an
algebraically closed field $\mathbb F$. One can associate to $G$ its
Lie algebra $L(G)$ and its algebra of distributions $U(G)$, which we
prefer to call the hyperalgebra of $G$. If $\mathbb F$ is of
characteristic zero, the hyperalgebra coincides with the universal
enveloping algebra $U(L(G))$ of $L(G)$, but this is not so in
positive characteristic.  $U(G)$ acts naturally on any $G$-module
and it turns out that, as conjectured originally by Verma and proved
by Sullivan \cite{sul}, every finite-dimensional $U(G)$-module can
be ``lifted'' to a rational finite-dimensional $G$-module. We will
restrict our attention to the case when $G$ is the Chevalley group
of adjoint type associated to a complex finite-dimensional simple
Lie algebra $\lie g$. In this case the algebra $U(G)$ is isomorphic
to the algebra $U(\lie g)_\mathbb F$ constructed by considering
Kostant's integral form of $U(\lie g)$ and tensoring with $\mathbb
F$ over $\mathbb Z$. It will suffice, for our purposes, to work over
the purely algebraic setting of $U(\lie g)_\mathbb F$.

Let $\lie g$ be as above and $\tlie g = \lie g\otimes\mathbb
C[t,t^{-1}]$ the loop algebra over $\lie g$. The finite-dimensional
representation theory of $\tlie g$ has been a very active research
topic in the last decades. It is related, for instance, to
integrable models and the Bethe ansatz in statistical mechanics. In
\cite{garala}, Garland introduced an integral form of $U(\tlie g)$
which can be used to construct what we call the hyper loop algebra
$U(\tlie g)_\mathbb F$ of $\lie g$ over $\mathbb F$ (see also
\cite{titsa,mitz}). The hyperalgebra $U(\lie g)_\mathbb F$ is naturally a
subalgebra of $U(\tlie g)_\mathbb F$.

The purpose of the present paper is to study some basic aspects of
the category $\tilde{\cal C}_\mathbb F$ of finite-dimensional
$U(\tlie g)_\mathbb F$-modules such as the classification of its
simple objects and its block decomposition when $\mathbb F$ is an
algebraically closed field of positive characteristic. In the case
$\mathbb F=\mathbb C$, thus $U(\tlie g)_\mathbb F = U(\tlie g)$,
these questions were studied in \cite{cint,cmsc,cpnew}. It turns out
that the simple finite-dimensional $\tlie g$-modules are
highest-weight modules with respect to the triangular decomposition
of $\tlie g$ obtained by ``looping'' the usual triangular
decomposition of $\lie g$. As usual,  we will call them
$\ell$-highest-weight modules to distinguish from those which are
highest-weight with respect to the triangular decomposition coming
from the Chevalley generators of $\tlie g$ (non-trivial
highest-weight representations with respect to the later
decomposition are infinite-dimensional). Moreover, all the simple
modules are isomorphic to suitable tensor products of the so-called
evaluation representations (obtained by pulling back the simple
$\lie g$-modules by the evaluation map $t\mapsto a$ for some nonzero
$a\in\mathbb C$). We prove that these two results hold in positive
characteristic, as well. This is done in Corollary \ref{ilgehw} and
Theorem \ref{tpd}, the later being a $U(\tlie g)_\mathbb F$-version
of Steinberg's Tensor Product Theorem. Using the tensor product
theorem we compute the dual representation of a given irreducible
one. For highest-weight representations with respect to the usual
triangular decomposition in positive characteristic see
\cite{garala,garalg,mat,mats,tits} and references therein.

The set of $\ell$-highest weights can be identified with ${\rm
rank}(\lie g)$-tuples of polynomials in $\mathbb F[u]$ with constant
term 1. For $\mathbb F=\mathbb C$, it was shown in \cite{cpweyl}
that there exists a family of universal finite-dimensional
$\ell$-highest-weight modules, called the Weyl modules. We prove
that the Weyl modules for $U(\tlie g)_\mathbb F$ can be defined in a
similar fashion when $\mathbb F$ is of positive characteristic. The
reason for calling these $\ell$-highest-weight modules
 Weyl modules comes from a conjecture in
\cite{cpweyl} stating that the Weyl modules for $U(\tlie g)$ can be
obtained as the classical limit of certain irreducible
finite-dimensional modules for the corresponding quantum loop
algebra, resembling the process of obtaining the Weyl modules for
$U(\lie g)_\mathbb F$ by reduction modulo $p$ of simple $\lie
g$-modules. This conjecture has been recently proved when $\lie g$
is of type $A$ in \cite{cl} using Gelfand-Tsetlin filtrations and
when $\lie g$ is of type $ADE$ in \cite{fl2} using Demazure modules.
Moreover, H. Nakajima has pointed out that the general case can be
deduced using the crystal and global basis results from
\cite{bn,kascb,kaslz,nakq,nake}. Other interesting related
references include \cite{cg,fl1,kmotu1,kmotu2}. We have an analogous
conjecture for the Weyl modules for $U(\tlie g)_\mathbb F$ when
$\mathbb F$ is of positive characteristic (Conjecture
\ref{cp=0}(a)), stating that they can be obtained from the Weyl
modules for $U(\tlie g)_{\mathbb F^0}$ by reduction modulo $p$,
where $\mathbb F^0$ is a suitable field of characteristic zero. As
$\mathbb Z$-lattices are easily seen not to be well-behaved with
respect to evaluation maps, we consider more general lattices for
this purpose. Namely, we consider lattices over the ring $\mathbb A$
of Witt vectors with coefficients in $\mathbb F$, after changing
scalars from $\mathbb C$ to the fraction field $\mathbb F^0$ of
$\mathbb A$.  We prove that all finite-dimensional
$\ell$-highest-weight $U(\tlie g)_{\mathbb F^0}$-modules whose coefficients
of the $\ell$-highest weights are in $\mathbb A$ and the leading ones are units in $\mathbb A$ contain an
admissible (ample) $\mathbb A$-lattice. Thus, we obtain all of
the irreducible modules as quotients of  modules coming from a
reduction modulo $p$ process. This is done in Theorem
\ref{latticeshg} and Corollary \ref{allirr}. Combining Conjecture
\ref{cp=0} with the one in \cite{cpweyl}, which is now a theorem as
remarked above, we have a bridge connecting the Weyl modules for
$U(\tlie g)_\mathbb F$ with certain irreducible representations for
quantum loop algebras (at generic quantization parameter).

As  corollaries of Conjecture \ref{cp=0}, we obtain a tensor product
decomposition of the Weyl modules  and the block decomposition of
$\widetilde{\cal C}_\mathbb F$. Although this tensor product
decomposition is the natural analogue of the one obtained in
\cite{cpweyl} for characteristic zero, the techniques used in that
paper do not seem to apply to our setting. In fact, our motivation
for considering the theory of reduction modulo $p$ originated from
the search for other methods which would lead to a proof of this
tensor product decomposition. Indeed we expect that this
decomposition holds in the context of $\mathbb A$-lattices
(Conjecture \ref{cp=0}(b)), thus allowing us to transfer the problem
to a characteristic zero setting. The block decomposition of
$\tilde{\cal C}_\mathbb F$ is described  similarly to that of
$\tilde{\cal C}_{\mathbb C}$ as well, i.e., the blocks are
parametrized by functions with finite support $\chi:\mathbb
F^\times\to P/Q$ called spectral characters. Here $P$ and $Q$ are
the weight and root lattices of $\lie g$, respectively, and $\mathbb
F^\times = \mathbb F-\{0\}$. The proof runs parallel to its
characteristic zero counterpart found in \cite{cmsc}, hence, the
tensor product decomposition for Weyl modules plays a key role.
However, our results on reduction modulo $p$ are needed in order to
both prove that the Weyl modules have a well-defined spectral
characters and obtain a positive characteristic version of
\cite[Proposition 3.4]{cmsc} -- a key ingredient in the construction
of certain useful indecomposable modules.

The paper is organized as follows. In section \ref{1stsec} we fix
the basic notation on finite-dimensional complex simple Lie algebras
and their loop algebras, define the hyperalgebras, and collect some
important Lemmas. Section \ref{U(g)Kir} is dedicated to a review of
the relevant facts about finite-dimensional $U(\lie g)_\mathbb
F$-modules. The main part of the paper consists of sections
\ref{mainsec} and \ref{redp}. In \ref{ellhw} we define
$\ell$-highest-weight modules and obtain the necessary relations
satisfied by the finite-dimensional ones. The classification of the
irreducible modules and the aforementioned tensor product and
duality results are done in \ref{s:cirhla}. The Weyl modules  are
constructed in \ref{wm}.  Section \ref{redp} ends the paper with the
results and the conjecture on reduction modulo $p$, as well as,
their  application to the description of the blocks.

\vskip25pt \noindent{\bf Acknowledgements:} D.J. is pleased to thank
the Max-Planck-Institut f\"ur Mathematik in Bonn for its hospitality
and support, FAPESP (processo 2006/00609-1) for supporting her visit
to the University of Campinas when part of this paper was developed,
and UNICAMP for its hospitality. The research of A.M. is partially
supported by CNPq (processo 303349/2005-0) and FAPESP (processo
2006/00833-9). We thank V. Chari for turning our attention to this
subject and for very helpful questions and suggestions. We also
thank P. Russel, A. Engler, and P. Brumatti for useful discussions
and pointers about discrete valuation rings and L. Scott for his
interest and helpful references.

\section{Hyperalgebras}\label{1stsec}

Throughout the paper  $\mathbb C, \mathbb Z,\mathbb Z_+,\mathbb N$ will denote the sets of complex  numbers, integers, non-negative integers, and positive integers, respectively. Given a ring $\mathbb A$, the underlying multiplicative group of units will be denoted by $\mathbb A^\times$. The dual of a vector space $V$ will be denoted by $V^*$.

\subsection{Preliminaries}

Let $I$ be the set of vertices of a finite-type connected Dynkin
diagram and let $\lie g$ be the associated simple complex Lie algebra with a fixed Cartan subalgebra $\lie h$ and nilpotent subalgebras $\lie n^{\pm}$. Denote by $R^+$ the set of positive roots so that
$$\lie n^\pm = \opl_{\alpha\in R^+}^{} \lie g_{\pm\alpha}, \quad\text{where}\quad \lie g_{\pm\alpha} = \{x\in\lie g: [h,x]=\pm\alpha(h)x, \ \forall \ h\in\lie h\}.$$
The simple roots will be denoted by $\alpha_i$, the fundamental weights by $\omega_i$, while $Q,P,Q^+,P^+$ will denote the root and weight lattices with corresponding positive cones, respectively. We equip $\lie h^*$ with the partial order $\lambda\le \mu$ iff $\mu-\lambda\in Q^+$. The Weyl group will be denoted by $\cal W$, its longest element by $w_0$, and the maximal positive root by $\theta$. Let $\langle\ ,\ \rangle$ be the bilinear form on $\lie h^*$ induced by the Killing form on $\lie g$ and, for $\lambda\in\lie h^*-\{0\}$, set $\lambda^\vee = 2\lambda/\langle\lambda,\lambda\rangle$ and $d_\lambda = \frac{1}{2}\langle\lambda,\lambda\rangle$. Then $\{\alpha_i^\vee:i\in I\}$ is a set of simple roots of the simple Lie algebra $\lie g^\vee$ whose Dynking diagram is obtained from that of $\lie g$ by reversing the arrows and $R^\vee=\{\alpha^\vee:\alpha\in R\}$ is its root system, where $R=R^+\cup (-R^+)$.  Moreover, if $\alpha = \sum_i m_i\alpha_i$ and $\alpha^\vee=\sum_i m_i^\vee\alpha_i^\vee$, then
\begin{equation}\label{m_ivee}
m_i^\vee = \frac{d_{\alpha_i}}{d_\alpha}m_i.
\end{equation}

If $\lie a$ is a Lie algebra (over any field $\mathbb F$), define its loop algebra $\tlie a=\lie a\otimes_{\mathbb F}  \mathbb F[t,t^{-1}]$ with bracket given by $[x \otimes t^r,y \otimes t^s]=[x,y] \otimes t^{r+s}$. Clearly $\lie a\otimes 1$ is a subalgebra of $\tlie a$ isomorphic to $\lie a$ and, by abuse of notation, we will continue denoting its elements by $x$ instead of $x\otimes 1$. In case $\lie a = \lie g$, we have $\tlie g = \tlie n^-\oplus \tlie h\oplus \tlie n^+$ and $\tlie h$ is an abelian subalgebra.

Let $U(\lie a)$ be the universal enveloping algebra of $\lie a$. Then $U(\lie a)$ is a subalgebra of $U(\tlie a)$ and, for $\lie a=\lie g$, multiplication establishes isomorphisms
$$U(\lie g)\cong U(\lie n^-)\otimes U(\lie h)\otimes U(\lie n^+) \qquad\text{and}\qquad U(\tlie g)\cong U(\tlie n^-)\otimes U(\tlie h)\otimes U(\tlie n^+).$$
The assignments $\triangle: \lie a\to U(\lie a)\otimes_\mathbb F
U(\lie a), x\mapsto x\otimes 1+1\otimes x$,  $S:\lie a\to \lie a,
x\mapsto -x$, and $\epsilon: \lie a\to \mathbb F, x\mapsto 0$, can
be uniquely extended so that $U(\lie a)$ becomes a Hopf algebra with
comultiplication $\triangle$, antipode $S$, and  counit $\epsilon$.
We shall denote by $U(\lie a)^0$ the augmentation ideal, i.e., the
kernel of $\epsilon$. Consider the associative $\mathbb F$-algebra
$U(\lie a)\otimes_\mathbb F\mathbb F[t,t^{-1}]$ with the obvious
tensor product structure and the usual bracket.  Clearly the
inclusion $\tlie a\hookrightarrow U(\lie a)\otimes_\mathbb F\mathbb
F[t,t^{-1}]$ is a Lie algebra map.
Therefore the next lemma is immediate from the
universal property of $U(\tlie a)$.

\begin{lem}\label{formalev}
There exists a unique algebra map $U(\tlie a)\to U(\lie a)\otimes_\mathbb F\mathbb F[t,t^{-1}]$ which is the identity on $\tlie a$.\hfill\qedsymbol
\end{lem}

We call the map given by this lemma the formal evaluation map and denote it by ${\rm ev}$. For each $a\in\mathbb F^\times$, consider the evaluation map $U(\lie a)\otimes_\mathbb F\mathbb F[t,t^{-1}]\to U(\lie a)$ sending $x\otimes f(t)$ to $f(a)x$ and denote by ${\rm ev}_a$ the composition of this map with ${\rm ev}$. Then ${\rm ev}_a$ is a surjective algebra homomorphism
$${\rm ev}_a:U(\tlie a)\to U(\lie a)$$
which we call the evaluation map at $a$.

\begin{rem}
Obviously, the existence of ${\rm ev}_a$ can  be proved similarly to the existence of ${\rm ev}$. However, we will use the formal evaluation map in order to prove the existence of evaluation maps in the context of hyper loop algebras using $\mathbb Z$-lattices only (cf. Proposition \ref{p:evmap} and the remark after Proposition \ref{p:homot}).
\end{rem}

\subsection{Reduction Modulo $p$} As usual, given any associative algebra $A$ over a field of characteristic zero, $a\in A$, and $k\in\mathbb Z_+$, we set $a^{(k)}=\frac{a^k}{k!}, \binom{a}{k} = \frac{a(a-1)\cdots(a-k+1)}{k!}\in A$.

Let $\Phi=\{ x^\pm_{\alpha}, h_{\alpha_i}: \alpha\in R^+,i\in I\}$
be a Chevalley basis for $\lie g$, where $x^\pm_\alpha\in \lie
g_{\pm\alpha}$, $h_\alpha=[x^+_\alpha,x^-_\alpha]$, and let
$x^\pm_{\alpha,r} = x^\pm_{\alpha}\otimes t^r, h_{\alpha,r} =
h_\alpha\otimes t^r$. When $r=0$ we may just write $x_{\alpha}^\pm$
and $h_{\alpha}$. Also, we may write
$x_{i,r}^\pm$ and $h_{i,r}$ in place of $x_{\alpha_i,r}^\pm$ and $h_{\alpha_i,r}$,
respectively. Notice that the set $\tilde\Phi = \{x_{\alpha,r}^\pm,
h_{i,r}:\alpha\in R^+,i\in I,r\in\mathbb Z\}$ is a basis for $\tlie
g$ and define $\tlie g_{\mathbb Z}$ to be the $\mathbb Z$-span of
$\tilde\Phi$. The $\mathbb Z$-span of $\Phi$ is a Lie $\mathbb
Z$-subalgebra of $\tlie g_{\mathbb Z}$ which we denote by $\lie
g_{\mathbb Z}$.

If $\mathbb F$ is any field, set
$$\lie g_{\mathbb F} = \lie g_{\mathbb Z}\otimes_{\mathbb Z} \mathbb F \qquad\text{and}\qquad \tlie g_{\mathbb F} = \tlie g_{\mathbb Z}\otimes_{\mathbb Z} \mathbb F,$$
so that $\lie g_{\mathbb F}$ and $\tlie g_{\mathbb F}$ are Lie algebras over $\mathbb F$.

Given $\alpha\in R^+, r\in\mathbb Z$, define elements $\Lambda_{\alpha,r}\in U(\tlie h)$ by the following equality of formal power series in $u$:
\begin{equation}
\Lambda_\alpha^\pm(u) = \sum_{r=0}^\infty \Lambda_{\alpha,\pm r} u^r = \exp\left( - \sum_{s=1}^\infty \frac{h_{\alpha,\pm s}}{s} u^s\right).
\end{equation}
We may write $\Lambda_{i,r}$ in place of $\Lambda_{\alpha_i,r}$. It follows from \eqref{m_ivee} that, if $\alpha=\sum_i m_i\alpha_i\in R^+$, then $h_\alpha = \sum_i m_i^\vee h_i$ and
\begin{equation}\label{Lambda_alpha}
\Lambda_\alpha^\pm(u) = \prod_{i\in I} (\Lambda_{\alpha_i}^\pm(u))^{m_i^\vee}.
\end{equation}
 We have (cf. \cite[Lemma 5.1]{garala}):
\begin{equation}\label{evLambda}
{\rm ev}(\Lambda_{\alpha,r}) = (-1)^r\binom{h_\alpha}{|r|}\otimes t^r.
\end{equation}
Set
\begin{equation}\label{H}
H_\alpha(u) = {\rm ev}_{-1}(\Lambda^+_\alpha(u)) = \exp\left( -\sum_{s\ge 1} h_\alpha \frac{(-u)^s}{s}\right),
\end{equation}
so that $(H_\alpha(u))_k$ (the coefficient of $u^k$ in $H_\alpha(u)$) is $\binom{h_\alpha}{k}$.

For $k\in\mathbb Z, k\ne 0$, consider also the endomorphism $\tau_k$
of $U(\tlie g)$ extending $t\mapsto t^k$ and set
$\Lambda_{\alpha,r;k}=\tau_k(\Lambda_{\alpha,r}),
\Lambda_{\alpha;k}^\pm(u) = \sum_{r=0}^\infty \Lambda_{\alpha,\pm
r;k} u^r $. Notice that $\binom{h_i}{k}$ is a polynomial in $h_i$ of
degree $k$. Hence,  the set $\{\binom{h_{1}}{k_1}\cdots \binom{h_{{\ell}}}{k_{\ell}}:k_j\in\mathbb Z_+\}$, where $\ell=|I|$, is a basis for $U(\lie h)$. Similarly, observe that $\Lambda_{i,\pm r;k},
r,k\in\mathbb N$,  is a polynomial in $h_{i,\pm sk}, 1\le s\le r$, whose leading term is $(-h_{i,\pm k})^{(r)}$. Finally, given an order on $\tilde \Phi$ and a PBW monomial with respect to this order,  we construct an ordered monomial  in the elements
$$(x_{\alpha,r}^\pm)^{(k)},\ \ \Lambda_{i,r;k},\ \ \binom{h_i}{k}, \qquad r,k\in\mathbb Z, k>0, \alpha\in R^+, i\in I,$$
 using the correspondence just discussed for the basis elements of $U(\tlie h)$, as well as, the obvious correspondence $(x_{\alpha,r}^\pm)^k \leftrightarrow  (x_{\alpha,r}^\pm)^{(k)}$. The set of monomials thus obtained is then a basis for $U(\tlie g)$, while the monomials involving only
$(x_{\alpha}^\pm)^{(k)},  \binom{h_i}{k}$ form a basis for
$U(\lie g)$. Let $U(\tlie g)_{\mathbb Z}$ (resp. $U(\lie g)_{\mathbb
Z}$) be the $\mathbb Z$-span of these monomials. The following
crucial theorem was proved in \cite{kosagz} ($U(\lie g)$ case) and
\cite{garala} ($U(\tlie g)$ case).

\begin{thm}\label{kostant}
$U(\tlie g)_{\mathbb Z}$ (resp. $U(\lie g)_{\mathbb Z}$) is the $\mathbb Z$-subalgebra of $U(\tlie g)$ generated by $\{(x_{\alpha,r}^\pm)^{(k)}, \alpha\in R^+,r,k\in\mathbb Z, k\ge 0\}$ (resp. $\{(x_{\alpha}^\pm)^{(k)}, \alpha\in R^+, k\in\mathbb Z_+\}$).\hfill\qedsymbol
\end{thm}

For $\lie a\in\{ \lie g, \lie n^\pm, \lie h, \tlie g,\tlie n^\pm,
\tlie h\}$, let $U(\lie a)_{\mathbb Z}$ denote the corresponding
$\mathbb Z$-subalgebra of $U(\tlie g)$. Given a field  $\mathbb F$, the $\mathbb F$-hyperalgebra of $\lie a$ is defined by
$$U(\lie a)_\mathbb F = U(\lie a)_{\mathbb Z}\otimes_{\mathbb Z} \mathbb F.$$
We will also refer to $U(\tlie g)_\mathbb F$ as the hyper loop algebra of $\lie g$ over $\mathbb F$. Then the PBW theorem gives
$$U(\lie g)_{\mathbb F}=U(\lie n^-)_{\mathbb F}U(\lie h)_{\mathbb F}U(\lie n^+)_{\mathbb F}\qquad\text{and}\qquad U(\tlie g)_{\mathbb F}=U(\tlie n^-)_{\mathbb F}U(\tlie h)_{\mathbb F}U(\tlie n^+)_{\mathbb F}.$$
Clearly, if $\mathbb F$ is of characteristic zero then $U(\lie a)_\mathbb F$ is naturally isomorphic to $U(\lie a_\mathbb F)$. For fields of positive characteristic we just have an algebra homomorphism $U(\lie a_\mathbb F)\to U(\lie a)_\mathbb F$ which is neither injective nor surjective.
If no confusion arises, we will write $x$, instead of $x\otimes 1$, for the image of an element $x\in U(\tlie g)_{\mathbb Z}$ in $U(\tlie g)_\mathbb F$.

Quite clearly the Hopf algebra structure on $U(\tlie g)$ preserves the $\mathbb Z$-forms $U(\lie a)_\mathbb Z$ and, therefore, induces a Hopf algebra structure on $U(\tlie g)_\mathbb F$ with counit given by $\epsilon((x_{\alpha,r}^\pm)^{(k)})=0, \epsilon(a)=a, a\in\mathbb F$, and comultiplication given by
\begin{equation}\label{comnil}
\triangle\left((x_{\alpha,r}^\pm)^{(k)}\right)= \sum_{l+m=k} (x_{\alpha,r}^\pm)^{(l)}\otimes (x_{\alpha,r}^\pm)^{(m)},
\end{equation}
\begin{equation}\label{comcart}
\triangle\left(\binom{h_i}{k}\right)= \sum_{l+m=k}\binom{h_i}{l}\otimes \binom{h_i}{m}, \quad\text{and}\quad
\triangle(\Lambda_{\alpha,\pm k}) = \sum_{l+m=k} \Lambda_{\alpha,\pm l}\otimes \Lambda_{\alpha,\pm m}.
\end{equation}
Moreover, the antipode on the basis of $U(\tlie h)_\mathbb F$ is determined by
\begin{equation}\label{antipode}
S(\Lambda^\pm_{\alpha;k}(u)) = (\Lambda^\pm_{\alpha;k}(u))^{-1} \qquad\text{and}\qquad S(H_\alpha(u)) = (H_\alpha(u))^{-1},
\end{equation}
where the inverses  in the last two equations are the ones of formal power series.

\subsection{Some Lemmas} We now collect some essential identities on $U(\tlie g)_\mathbb F$, when $\mathbb F$ is a field of characteristic $p>0$.
We begin with the following trivial observation:
\begin{equation}\label{binid}
(x_{\alpha,r}^\pm)^{(k)}(x_{\alpha,r}^\pm)^{(l)} = \binom{k+l}{k}(x_{\alpha,r}^\pm)^{(k+l)}.
\end{equation}
From this, one easily deduces
\begin{equation}\label{ppowerx}
((x_{\alpha,r}^\pm)^{(k)})^p=0.
\end{equation}
It is well-known (see \cite{humhyp})  that the elements $\binom{h_i}{k}$ satisfy
\begin{equation}\label{ppowerh}
\binom{h_i}{k}^p = \binom{h_i}{k}
\end{equation}
and it is easy to see that we have
\begin{equation}\label{comutxh}
\binom{h_i}{l}(x_{\alpha,r}^{\pm})^{(k)} = (x_{\alpha,r}^{\pm})^{(k)}\binom{h_i\pm k\alpha(h_i)}{l}.
\end{equation}

Given $\alpha\in R^+, s\in\mathbb Z$, define
$$X^-_{\alpha;s,\pm}(u) = \sum_{r\ge 1} x^-_{\alpha,\pm(r+s)}u^r.$$

\begin{lem}
We have:
\begin{gather}\label{koslem}
(x_{\alpha}^+)^{(l)}(x_{\alpha}^-)^{(k)} = \sum_{m=0}^{\min\{k,l\}} (x_{\alpha}^-)^{(k-m)}\binom{h_\alpha-k-l+2m}{m}(x_{\alpha}^+)^{(l-m)}\\
\notag\text{and}\\ \label{basicrel}
(x^+_{\alpha,\mp s})^{(l)}(x^-_{\alpha,\pm(s+1)})^{(k)} \in (-1)^l \left((X_{\alpha;s,\pm}^-(u))^{(k-l)}\Lambda_{\alpha}^{\pm}(u)\right)_k + U(\tlie g)_\mathbb F U(\tlie n^+)_\mathbb F^0.
\end{gather}
In \eqref{basicrel}, $k\ge l\ge 1$ and the subindex $k$  means the coefficient of $u^k$ of the above power series.
\end{lem}

\begin{proof}
It suffices to prove that the relations hold in $U(\tlie g)_\mathbb Z$. For both claims the strategy is to commute the elements on the left hand side. The proof of \eqref{koslem} can be found in  \cite[Lemma 26.2]{humb}.

Relation \eqref{basicrel} was proved in \cite[Lemma 7.5]{garala} for $s=0$ and the choice of ``$\pm$'' such that we have ``+'' on the right hand side. Consider the subalgebra of $U(\tlie g)_\mathbb Z$ generated by $(x_{\alpha,r}^\pm)^{(k)}$ for a fixed $\alpha\in R^+$. It is easy to see that, for each $s\in\mathbb Z$, the assignment $(x_{\alpha,r}^\pm)^{(k)}\mapsto (x_{\alpha,r\pm s}^\pm)^{(k)}$ extends uniquely to an algebra automorphism of this subalgebra which is the identity when restricted to $U(\tlie h)_\mathbb Z$. The general case of \eqref{basicrel} (with ``$+$'' on the right hand side) follows easily from the case $s=0$ using these automorphisms (see also \cite[Lemma 1.3]{cpweyl}). For the opposite choice of ``$\pm$'', just apply the automorphism determined by the assignment $(x_{\alpha,r}^\pm)^{(k)}\mapsto (x_{\alpha,-r}^\pm)^{(k)}$.
\end{proof}

The following lemmas will be needed in the proof of Theorem
\ref{wmfd}. Consider monomials involving only the elements
$(x_{\alpha,r}^-)^{(k)}$. Define the degree of
$(x_{\alpha,r}^-)^{(k)}$ to be $k$ and extend it additively.

\begin{lem}\label{sbrw}
Let $\alpha,\beta\in R^+, k,l\in\mathbb Z_+, r,s\in\mathbb Z$.  Then $(x_{\alpha,r}^-)^{(k)}(x_{\beta,s}^-)^{(l)}$ is in the span of $(x_{\beta,s}^-)^{(l)}(x_{\alpha,r}^-)^{(k)}$ together with monomials of degree strictly smaller than $k+l$.
\end{lem}

\begin{proof}
Immediate from the proof of \cite[Lemma 26.3.C]{humb} using that
$U(\tlie n^-)_{\mathbb Z}$ is $(Q^+\times\mathbb Z)$-graded.
\end{proof}

The next lemma is part of \cite[Lemma 5.11]{garala} and shows that the elements $\Lambda_{\alpha,r;k}$ are linear combinations of products of the elements $\Lambda_{\alpha,s}$.

\begin{lem}\label{ht^snot0}
In $U(\tlie h)_{\mathbb Z}$, for all $k,s\in \mathbb N$ and $\alpha\in R^+$, we have
\begin{gather*}
\Lambda_{\alpha,\pm s;k} = k\Lambda_{\alpha,\pm sk}+\sum_{(\vec r,\vec n)} m_{\vec r,\vec n}\Lambda_{\alpha,\pm r_1}^{n_1}\cdots\Lambda_{\alpha,\pm r_l}^{n_l},
\end{gather*}
for some $m_{\vec r,\vec n}\in\mathbb Z$. The sum is over the pairs $(\vec r,\vec n)$ where $\vec r = (r_1,\cdots, r_l)$ and $\vec n=(n_1,\cdots, n_l)$ are such that $l,r_j,n_j\in \mathbb N$,  $r_i\ne r_j, l\sum_j n_j>1$, and $\sum_j n_jr_j=sk$.\hfill\qedsymbol
\end{lem}

\subsection{Frobenius Homomorphism}\label{frsec}

\begin{lem}\label{lessgen}
$U(\tlie g)_\mathbb F$ (resp. $U(\lie g)_\mathbb F$) is generated as an algebra by the elements $(x_{\alpha,r}^\pm)^{(p^k)}$ (resp. $(x_{\alpha}^\pm)^{(p^k)}$), $\alpha\in R^+, k,r\in\mathbb Z, k\ge 0$. Moreover, $U(\lie h)_\mathbb F$ is generated as an algebra by $\binom{h_i}{p^k}, i\in I, k\in\mathbb Z_+$.
\end{lem}

\begin{proof}
The first statement  is immediate from Theorem \ref{kostant} and \eqref{binid}. The second is a statement on $U(\lie g)_\mathbb F$ and is well known.
\end{proof}

It is known that there exists a Hopf algebra map $\tilde\phi: U(\tlie g)_\mathbb F\to U(\tlie g)_\mathbb F$ sending  $(x_{\alpha,r}^\pm)^{(p^k)}$ to  $(x_{\alpha,r}^\pm)^{(p^{k-1})}$ (with the convention that the later is zero when $k=0$). We will denote the restriction of $\tilde\phi$ to $U(\lie g)_\mathbb F$ by $\phi$ and call both of them the (arithmetic) Frobenius homomorphisms. The first formula below is well known and the second was proved in \cite{cpswitt}.
\begin{equation}\label{freq}
\phi\left(\binom{h_i}{p^k}\right) = \binom{h_i}{p^{k-1}} \qquad\text{and}\qquad \tilde\phi(\Lambda_{i,r}) =
\begin{cases}
\Lambda_{i,r/p}, & \ \text{if } p\text{ divides }r\\
0, & \text{otherwise}.
\end{cases}
\end{equation}

The proof of the existence of the map $\phi$ can be found in
\cite{janb}, for instance. For the existence of $\tilde\phi$, see
\cite[Lemma 1.3]{mats} and \cite[Lemma 9.5]{cpr}.

Given a $U(\tlie g)_\mathbb F$-module (resp.  $U(\lie g)_\mathbb F$-module) $V$, we denote
by $V^{\tilde\phi^m}$ (resp. $V^{\phi^m}$) the pull-back of $V$ by
${\tilde\phi}^m$ (resp. $\phi^m$).

\section{Review of Finite-Dimensional $U(\lie g)_\mathbb F$-Modules}\label{U(g)Kir}

In this section we review some results on finite-dimensional
representations of $U(\lie g)_\mathbb F$ which will be relevant for our
purposes. In the first subsection we consider the case $\mathbb F=\mathbb
C$, where we summarize the basic results without proofs. The
literature for this subsection is vast and well known (all the
results we mention can be found in \cite{humb} to name but one
reference). In the other subsections $\mathbb F$ will be an algebraically
closed field of characteristic $p>0$. Essentially all of the results
can be found in \cite{janb} (see also \cite{lnm131}), although the
approach there is heavily geometric. Our approach follows  that of
\cite[Chapter VII]{humb} and \cite{humhyp}. Since some proofs are
relevant for  section \ref{mainsec}, we consider it appropriate to
sketch them.

\subsection{Characteristic Zero and Lattices} Given a $U(\lie g)$-module $V$, a vector $v\in V$ is called a weight vector if $hv=\mu(h)v$ for some $\mu\in\lie h^*$ and all $h\in\lie h$. The subspace consisting of weight vectors of weight $\mu$ will be denoted by $V_\mu$. If $v$ is a weight vector such that $\lie n^+v = 0$, then $v$ is called a highest-weight vector. If $V$ is generated by a highest-weight vector of weight $\lambda$, then $V$ is said to be a highest-weight module of highest weight $\lambda$.

The following theorem summarizes the basic facts about
finite-dimensional $U(\lie g)$-modules.

\begin{thm}\label{c0cig} Let $V$ be a finite-dimensional $U(\lie g)$-module.
\begin{enumerate}
\item $V=\opl_{\mu\in\lie h^*}^{} V_\mu$ and $\dim V_\mu = \dim V_{w\mu}$ for all $w\in\cal W$.
\item $V$ is completely reducible.
\item For each $\lambda\in P^+$ the $U(\lie g)$-module $V^0(\lambda)$ generated by a vector $v$ satisfying
$$x_i^+v=0, \quad h_iv=\lambda(h_i)v, \quad (x_i^-)^{\lambda(h_i)+1}v=0,\quad\forall\ i\in I,$$
is irreducible and finite-dimensional. If $V$ is irreducible, then
$V$ is isomorphic to $V^0(\lambda)$ for some $\lambda\in P^+$.
\item If $\lambda\in P^+$ and $V\cong V^0(\lambda)$, then $V_\mu\ne 0$ iff $w\mu\le\lambda$ for all $w\in\cal W$. Furthermore, the minimal weight of $V^0(\lambda)$ is $w_0\lambda$.
\hfill\qedsymbol
\end{enumerate}
\end{thm}

An admissible lattice for a $U(\lie g)$-module $V$ is the $\mathbb
Z$-span of a basis for $V$ which is invariant under the action of
$U(\lie g)_{\mathbb Z}$.  The basic results about lattices can be
summarized in the following Theorem (see \cite{humb}).

\begin{thm}\label{elat} Let $V,W$ be  finite-dimensional $U(\lie g)$-modules.
\begin{enumerate}
\item If $L$ is an additive subgroup of $V$ which is invariant under the action of $U(\lie g)_{\mathbb Z}$, then $L=\opl_{\mu\in P}^{} L_\mu$, where $L_\mu= L\cap V_\mu$.
\item There exists an admissible lattice for $V$.
\item If $L,M$ are admissible lattices for $V,W$, respectively, then $L\otimes_{\mathbb Z} M$ is an admissible lattice for $V\otimes W$.
\item If $V$ is an irreducible module and $v$ is a highest-weight vector of weight $\lambda$, then $L=U(\lie n^-)_{\mathbb Z}v$ is minimal in the set of admissible  lattices for $V$ satisfying $L_\lambda = \mathbb Zv$. \hfill\qedsymbol
\end{enumerate}
\end{thm}

\subsection{Classification of Irreducible Modules in Positive Characteristic} From now on, $\mathbb F$ is an algebraically closed field of characteristic $p>0$ and $\mathbb F_p$ denotes its prime field. In the present subsection we recall the methods used to classify the irreducible representations of $U(\lie g)_\mathbb F$ up to isomorphism. Although the classification is the same as in the case of $U(\lie g)$, the methods are quite different and will be used when we treat the case of $U(\tlie g)_\mathbb F$.

Let $V$ be a $U(\lie h)_\mathbb F$-module. A nonzero vector $v\in V$ is
called a weight vector if there exists $\gb z = (z_{i,k}),
z_{i,k}\in \mathbb F, i\in I,k\in\mathbb Z_+$, called the weight of $v$, such that  $\binom{h_i}{p^k}v
= z_{i,k}v$. Notice that \eqref{ppowerh} implies that $z_{i,k}$ must be in $\mathbb F_p$. We say
that $\gb z$ is  integral (resp. dominant integral) if $z_{i,k} =
\binom{\mu(h_i)}{p^k}$ for some $\mu\in P$ (resp. $\mu\in P^+$). In
that case we identify $\gb z$ with $\mu$ and say that $v$ has weight
$\mu$.
If $V$ is a $U(\lie g)_\mathbb F$-module and $v$ is weight vector such that
$(x_\alpha^+)^{(k)}v=0$ for all $\alpha\in R^+, k\in\mathbb N$, then $v$ is
said to be  a highest-weight  vector. If $V$ is generated  by a
highest-weight vector, $V$ is called a highest-weight module.

Since the $\binom{h_i}{p^k}$ commute, we can decompose any
finite-dimensional representation $V$ of $U(\lie g)_\mathbb F$ in a direct
sum of generalized eigenspaces for the action of $U(\lie h)_\mathbb F$:
$$V = \opl_{\gb z}^{} V_{\gb z}.$$
 We say that $\gb z$ is a weight of $V$ if $V_{\gb z}\ne 0$ and, in that case, $V_{\gb z}$ is called a weight space of $V$. In the case $\gb z$ is integral we write $V_\mu$ instead of $V_{\gb z}$.

Given $\gb z = (z_{i,k})$ and $\mu\in P$ define $\gb z+\mu=\gb y$ by the equality $y_{i,k}v=\binom{h_i+\mu(h_i)}{p^k}v$, where $v$ is some weight vector of weight $\gb z$. It follows from \eqref{comutxh} that if $v$ has weight $\gb z$ then
$(x_\alpha^\pm)^{(k)}v$ is either zero or has weight $\gb z\pm
k\alpha$. Hence, if $v$ is a highest-weight vector for a
highest-weight representation $V$, we have $\dim(V_{\gb z})=1$ and
$V_{\gb y}\ne 0$ only if $\gb y\le \gb z$, where $\gb y\le \gb z$
iff $\gb y = \gb z-\eta$ for  some $\eta\in Q^+$. Standard arguments then
show:

\begin{prop}\label{p:hq}
Every highest-weight module is indecomposable and has a unique
maximal proper submodule, hence, also a unique irreducible
quotient.\hfill\qedsymbol
\end{prop}

Recall that any nonnegative integer $m$ can be written uniquely as $m = \sum_{j\ge 0} m_jp^j$, where $0\le m_j<p$, so that $\binom{m}{p^r} = m_r \ ({\rm mod}\ p)$ for all $r\ge 0$. We shall write $\overline m$ for the image of $m\in \mathbb Z$ in $\mathbb F_p$.

\begin{thm}\label{cig}
If $V$ is an irreducible finite-dimensional $U(\lie g)_\mathbb F$-module,
then $V$ is a highest-weight representation with dominant integral
highest weight.
\end{thm}

\begin{proof} Since $V$ is irreducible, the generalized eigenspaces spaces $V_{\gb z}$ are in fact eigenspaces. Moreover, since $V$ is finite-dimensional, it also follows that there exists a maximal weight $\gb z$ and, hence, $V$ is a highest-weight module. It remains to prove that $\gb z$ is dominant integral. Let $v$ be a highest-weight vector for $V$. As we have already observed above, $(x_\alpha^-)^{(k)}v$ is either zero or has weight $\gb z-k\alpha$. This implies that, for every $i\in I$, there exists $N_i\in\mathbb Z_+$ minimal such that $(x_i^-)^{(p^k)}v = 0$ for all $k\ge N_i$. Moreover, we conclude from \eqref{koslem} with $k=l\ge p^{N_i}$ that $\binom{h_i}{p^{r}}v=0$ for all $r\ge N_i$. Now we easily see that $\gb z$ coincides with $\lambda\in P^+$ defined by $\lambda(h_i) = \sum_{j=0}^{N_i-1} m_{i,j}p^j$ with $0\le m_{i,j}<p$ such that $\overline m_{i,r}=z_{i,r}$.
\end{proof}

In order to complete the classification of the irreducible $U(\lie
g)_\mathbb F$-modules in terms of highest weights, it remains to prove
that for every $\lambda\in P^+$, there exists an irreducible $U(\lie
g)_\mathbb F$-module having $\lambda$ as highest weight. We will use reduction modulo $p$ as follows.

Let $V$ be  a finite-dimensional $U(\lie g)$-module and $L$ an admissible lattice for $V$.  Setting $L_\mathbb F = L\otimes_{\mathbb Z} \mathbb F$, we have  that $L_\mathbb F$ is a $U(\lie g)_\mathbb F$-module and $\dim_\mathbb F(L_\mathbb F) = \dim_{\mathbb C}
(V)$. The $U(\lie g)_\mathbb F$-module $L_\mathbb F$ is called a reduction modulo $p$ of $V$ (via $L$). If $L$ is a minimal admissible lattice for $V=V^0(\lambda)$, then
$L_\mathbb F$ is clearly highest-weight with highest weight $\lambda$. Hence, by Proposition \ref{p:hq}, it has a finite-dimensional irreducible quotient. Let $V(\lambda)$ denote this quotient.

We end this subsection remarking that the following statement remains true in positive characteristic.

\begin{prop}\label{wgi}
Let $V$ be a finite-dimensional $U(\lie g)_\mathbb F$-module. The
generalized eigenspaces $V_\mu$ are in fact eigenspaces and $\dim
V_\mu = \dim V_{w\mu}$ for all $w\in\cal W$.\hfill\qedsymbol
\end{prop}

\subsection{Weyl Modules and Duality}

\begin{defn}
Given $\lambda\in P^+$,  let $W(\lambda)$ be the $U(\lie
g)_\mathbb F$-module generated by a vector $v$ satisfying
\begin{gather}
(x_{\alpha}^+)^{(p^k)}v = 0, \quad \binom{h_i}{p^k}v = \binom{\lambda(h_i)}{p^k}v,\quad (x_{\alpha}^-)^{(l)}v = 0, \quad \forall\ \alpha\in R^+, i\in I, k,l\in\mathbb Z_+, l>\lambda(h_\alpha).
\end{gather}
\end{defn}

The modules $W(\lambda)$ are called Weyl modules. One can show that
every finite-dimensional highest-weight $U(\lie g)_\mathbb F$-module
is a quotient of some $W(\lambda)$.   A comparison between the
definition of $W(\lambda)$ and Theorem \ref{c0cig}(c) hints that we
have the following theorem which is a consequence of Kempf's
Vanishing Theorem and shows that $W(\lambda)$ is universal in the
family of finite-dimensional highest-weight modules with highest
weight $\lambda$.

\begin{thm}\label{wmgfd}
Let $\lambda\in P^+$ and $L$ a minimal admissible lattice for $V^0(\lambda)$. Then $W(\lambda)$ is isomorphic to $L_\mathbb F$. In particular, $W(\lambda)$ is finite-dimensional.\hfill\qedsymbol
\end{thm}

The notions of lowest-weight vector and lowest-weight module are
defined similarly to the corresponding highest-weight notions. It is
well-known that $V^0(\lambda)$ is a lowest-weight module with lowest
weight $w_0\lambda$, where $w_0$ is the longest element of $\cal W$.
Given a highest-weight vector $v$ for $V^0(\lambda)$ and a reduced
expression $w_0=s_{i_l}\cdots s_{i_1}$,  set $m_{i_k}\in\mathbb Z_+,
k=1,\dots, l$, to be $(s_{i_{k-1}}\cdots s_{i_1}\lambda)(h_{i_k})$.
Then a lowest-weight vector of $V^0(\lambda)$ is given by
$v'=(x_{i_l}^-)^{(m_{i_l})}\cdots (x_{i_1}^-)^{(m_{i_1})}v$ and,
moreover,  $v=(x_{i_1}^+)^{(m_{i_1})}\cdots
(x_{i_l}^+)^{(m_{i_l})}v'$. This shows that the image of $v'$ in the
irreducible quotient of $W(\lambda)$ is nonzero and moreover:

\begin{cor}\label{loww}
For all $\lambda\in P^+$, $W(\lambda)$ and  $V(\lambda)$ are
lowest-weight modules with lowest weight
$w_0\lambda$.\hfill\qedsymbol
\end{cor}

Since $U(\lie g)_{\mathbb F}$ is a Hopf algebra, given a $U(\lie
g)_\mathbb F$-module $V$, one can equip the dual vector space $V^*$
with a structure of $U(\lie g)_\mathbb F$-module where the action of
$x\in U(\lie g)_\mathbb F$ on $f\in V^*$ is given by
$$(xf)(v) = f(S(x)v)$$
for all $v\in V$.

\begin{prop}\label{p:dualfin}
Let $V$ be a finite-dimensional $U(\lie g)_\mathbb F$-module. Then
\begin{enumerate}
\item The natural isomorphism of vector spaces $V\to V^{**}$ is a $U(\lie g)_\mathbb F$-module isomorphism.
\item If $V=V(\lambda),\lambda\in P^+$, then $V^*\cong V(-w_0\lambda)$.
\end{enumerate}
\end{prop}

\begin{proof}
Part (a) is an immediate consequence of the fact that $S^2$ is the identity. Now if $V$ is irreducible, since duality preserves exact sequences, it follows from (a) that $V^*$ is also irreducible. From \eqref{antipode} we conclude that $V_\mu\ne 0$ iff $V^*_{-\mu}\ne 0$. The final claim  is now immediate from Corollary \ref{loww}.
\end{proof}

\subsection{Tensor Products}

We now recall Steinberg's Tensor Product Theorem
\cite{sterag}. We sketch only the part of the proof which will be relevant for section \ref{mainsec}. Our argument  essentially follows the one  given in \cite{cps}.  Let $P_p^+=\{\lambda\in P^+: \lambda(h_i)<p, \forall\
i\in I\}$.  We shall use the following lemma and refer to the aforementioned references for its proof.

\begin{lem}\label{redtp}
Let $\lambda,\mu\in P^+_p-\{0\}$. Then $V(\lambda)$ is irreducible
as $\lie g_\mathbb F$-module and $V(\lambda)\otimes V(\mu)$ is reducible as
$U(\lie g)_\mathbb F$-module.\hfill\qedsymbol
\end{lem}

\begin{thm}\label{steinberg}
For $\lambda\in P^+$, let $\lambda_k$ be the unique elements of $P^+_p$ such that $\lambda = \sum_{k=0}^m p^k\lambda_k$. Then
$V(\lambda) \cong \otimes_k V(p^k\lambda_k)$. Moreover, if $\mu_j\in P^+_p-\{0\}$ and $l_j\in\mathbb Z_+, j=0,\cdots, n$, are such that $\otimes_{j=0}^{n} V(p^{l_j}\mu_j)\cong V(\lambda)$, then $m=n$ and (up to reordering) $\mu_k=\lambda_k$ and $l_k=k$ for all $k$.
\end{thm}

\begin{proof}
First observe that for any $\mu\in P^+_p$ and $k\in\mathbb Z_+$  we
have $V(p^k\mu)\cong V(\mu)^{\phi^k}$ (see section \ref{frsec}).
Therefore, $(x^\pm_{\alpha})^{(p^l)} V(p^k\mu) = 0$ if $l< k$. Now
let $v_k$ be highest-weight vectors for $V(p^k\lambda_k)$,
$V'=\otimes_{k=1}^m V(p^k\lambda_k)$, and $v=\sum_{i} w_i\otimes
w'_i\in V(\lambda_0)\otimes V'$, where $w'_i$ are linearly
independent. Then $x_{\alpha}^+v = \sum_{i} (x_{\alpha}^+w_i)\otimes
w'_i$. Since $V(\lambda_0)$ is irreducible as $\lie g_\mathbb F$-module, it
follows that $x_\alpha^+v=0$ only if $v = v_0\otimes v'$ for some
$v'\in V'$. Now let $V''=\otimes_{k=2}^m V(p^k\lambda_k)$, and
$v=v_0\otimes (\sum_{i} w'_i\otimes w''_i)\in V(\lambda_0)\otimes
V(p\lambda_1)\otimes V''$, where $w''_i$ are linearly independent.
Then $(x_{\alpha}^+)^{(p)}v = v_0\otimes(\sum_{i}
((x_{\alpha}^+)^{(p)}w'_i)\otimes w''_i)$. Since $V(\lambda_1)$ is
irreducible as $\lie g_\mathbb F$-module, it follows that
$(x_\alpha^+)^{(p)}v=0$ only if $v = v_0\otimes v_1\otimes v''$ for
some $v''\in V''$. Continuing like this we see that $\otimes_{k=0}^m
V(p^k\lambda_k)$ is irreducible. Since it is clearly a
highest-weight module with highest-weight $\lambda$, the first
statement is proved. On the other hand, we must  have $\lambda_k=
\sum_{j\in J_k} \mu_j$ where $J_k=\{j:l_{j}=k\}$. Therefore, if
$\{\mu_j\}$ were not as stated, there would clearly exist $j\ne j'$
such that $l_{j}=l_{j'}$. The lemma above would then imply
$V(p^{l_j}\mu_j)\otimes V(p^{l_{j'}}\mu_{j'})$ is reducible and,
hence, also  $\otimes_{j=0}^{n} V(p^{l_j}\mu_j)$.
\end{proof}

\begin{rem}
One of the reasons Theorem \ref{steinberg} is important comes from
the fact that the (finitely many) modules $V(\lambda), \lambda\in P_p^+$, are irreducible as modules
for the subalgebra of $U(\lie g)_\mathbb F$ generated by $x_{\alpha}^{\pm}$
(since they are irreducible as $\lie g_\mathbb F$-modules). This is a finite-dimensional algebra, called the restricted universal enveloping algebra of $\lie g_\mathbb F$.
Hence, the study of finite-dimensional irreducible $U(\lie g)_\mathbb F$-modules is reduced to the study of finitely many modules for a finite-dimensional algebra.
\end{rem}

\section{Finite-Dimensional $U(\tlie g)_\mathbb F$-Modules}\label{mainsec}

In this section we establish some basic results about the category of finite-dimensional $U(\tlie g)_\mathbb F$-modules such as the classification of the irreducible ones and the characterization of the universal highest-weight modules.

\subsection{$\ell$-Highest-Weight Modules}\label{ellhw}

Let $V$ be a $U(\tlie g)_\mathbb F$-module. We say $v\in V$ is an
$\ell$-weight vector if it is an eigenvector for the action of
$U(\tlie h)_\mathbb F$, i.e., if there exist $z_{i,k},\varpi_{i,r}\in \mathbb F$
such that
\begin{gather}
\binom{h_i}{p^k}v = z_{i,k}v, \qquad \Lambda_{i,r}v = \varpi_{i,r} v,
\end{gather}
for all  $i\in I$ and all $r,k\in\mathbb Z, k\ge 0$. In that case
the corresponding functional $\gb\varpi\in (U(\tlie h)_\mathbb F)^*$ is
called the $\ell$-weight of $v$. If $v$ is an $\ell$-weight vector
and $(x_{\alpha,r}^+)^{(k)}v = 0$ for all $\alpha\in R^+$ and all
$r,k\in\mathbb Z, k>0$, we say $v$ is an $\ell$-highest-weight
vector. If $V$ is generated by an $\ell$-highest-weight vector, we
say $V$ is an $\ell$-highest-weight module.

Given a  finite-dimensional $U(\tlie g)_\mathbb F$-module $V$ we know from
section \ref{U(g)Kir} that $V$ can be written as the direct sum of
its weight spaces when regarded as $U(\lie g)_\mathbb F$-module:
$$V = \opl_{\mu\in P}^{} V_\mu.$$
 Moreover, since $U(\tlie h)_\mathbb F$ is a commutative algebra, we can also write the following decomposition of $V$ into direct sum of generalized eigenspaces for the action of $U(\tlie h)_\mathbb F$:
$$V = \opl_{\gb\varpi\in (U(\tlie h)_\mathbb F)^*}^{} V_{\gb\varpi}.$$

The next proposition establishes a set of relations  satisfied by all finite-dimensional
$\ell$-highest-weight modules.

\begin{prop}\label{ellhwrel}
Let $V$ be a finite-dimensional $U(\tlie g)_\mathbb F$-module, $\lambda\in
P^+$,  and $v\in V_\lambda$ be such that
$$(x_{\alpha,s}^+)^{(k)}v=0 \qquad\text{and}\qquad \Lambda_{i,s}v = \omega_{i,s}v,$$
for all $\alpha\in R^+,i\in I, k,s\in\mathbb Z, k>0$, and some $\omega_{i,s}\in \mathbb F$. Then
$$(x_{\alpha,s}^-)^{(k)}v = \Lambda_{i,\pm r}v=0 \quad \text{for all } k>\lambda(h_\alpha), r>\lambda(h_i), s\in\mathbb Z.$$ Moreover, $\omega_{i,\pm\lambda(h_i)}\ne 0$ and there exist polynomials $f_i\in\mathbb F[t_0,t_1,\cdots,t_{\lambda(h_i)}]$, depending only on $\lambda(h_i)$, such that
$$\omega_{i,-r} = f_i(\omega_{i,\lambda(h_i)}^{-1}, \omega_{i,1},\cdots, \omega_{i,\lambda(h_i)})$$
 for all $r=1,\cdots, \lambda(h_i)$.
\end{prop}

\begin{proof}
For each $r\in\mathbb Z, \alpha\in R^+$, the elements
$(x_{\alpha,\pm r}^{\pm})^{(k)},  k\in\mathbb Z_+$, generate a
subalgebra $U(\tlie g_{\alpha,r})_\mathbb F$ of $U(\tlie g)_\mathbb F$ isomorphic to $U(\lie{sl}_2)_\mathbb F$. Hence,
the equality $(x_{\alpha,r}^-)^{(k)}v = 0$ for $k>\lambda(h_\alpha)$
follows from the fact that $v$ generates a (finite-dimensional)
highest-weight module for this subalgebra, which is then isomorphic
to a quotient of the Weyl module $W(\lambda(h_\alpha))$.

Setting $\alpha=\alpha_i,s=0, l=k=r$ in \eqref{basicrel}  we get $\Lambda_{i,\pm r}v=0$ for $r>|\lambda(h_i)|$. Now, choosing $r=\lambda(h_i)$, we see that $\omega_{i,\pm\lambda(h_i)}\ne 0$ . In fact, since $W=U(\tlie g)_\mathbb Fv$ is a finite-dimensional representation for $U(\lie g)_\mathbb F$ having $W_\lambda=\mathbb Fv$ as its highest-weight space by equation \eqref{comutxh}, it follows that $\lambda-(r+m)\alpha_i$ is not a weight of $W$ for any $m>0$. Therefore $(x^-_i)^{(m)}(x_{i,\pm 1}^-)^{(r)}v=0$ for all $m\in\mathbb N$. On the other hand, by considering the subalgebra $U(\tlie g_{\alpha_i,\mp 1})_\mathbb F$, we see that $(x_{i,\pm 1}^-)^{(r)}v\ne 0$. It follows that $(x_{i,\pm 1}^-)^{(r)}v$ generates a lowest weight finite-dimensional representation of $U(\tlie g_{\alpha_i,0})_\mathbb F$ and, in particular, $0\ne (x_{i}^+)^{(r)}(x_{i,\pm 1}^-)^{(r)}v=\Lambda_{i,\pm r}v$.

For the last statement we proceed by induction on $r=1,\cdots,\lambda(h_i)=m_i$.  Setting $\alpha=\alpha_i, s=0,l=m_i$ and $k=l+r$ in \eqref{basicrel} we get
$$\omega_{i,m_i}(x_{i,1}^-)^{(r)}v + \sum_{j=1}^{m_i} \omega_{i,m_i-j}Y_jv=0,$$
where $\omega_{i,0}=1$ and $Y_j$ is the sum of the monomials $(x_{i,1}^-)^{(k_1)}\cdots (x_{i,m_i+1}^-)^{(k_{m_i+1})}$ such that $\sum_n k_n = r$ and $\sum_n nk_n = r+j$. Now, since $-r<r+j-2r<m_i$, it is not difficult to see that $(x_{i,-2}^+)^{(r)}Y_j \in U(\tlie g)_\mathbb FU(\tlie n^+)_\mathbb F^0 + H_j$, where $H_j$ is a linear combination of monomials of the form $\Lambda_{i,r_1}\cdots \Lambda_{i,r_m}$ such that $-r<r_j< m_i$. Moreover, $(x_{i,-2}^+)^{(r)}(x_{i,1}^-)^{(r)} \in (-1)^r\Lambda_{i,-r} + U(\tlie g)_\mathbb FU(\tlie n^+)_\mathbb F^0$ by \eqref{basicrel}. Hence,
$$0=(x_{i,-2}^+)^{(r)}\left(\omega_{i,m_i}(x_{i,1}^-)^{(r)}v + \sum_{j=1}^{m_i} \omega_{i,m_i-j}Y_jv\right)= (-1)^s\omega_{i,-r}\omega_{i,m_i}v + \sum_{j=1}^{m_i} \omega_{i,m_i-j}H_jv.$$
From here it is easy to deduce the last statement.
\end{proof}

We would like to be more precise about the last statement of the
previous proposition (cf. \cite[Proposition 1.1(v)]{cpweyl}).
Namely, we want to prove that
\begin{equation}\label{e:Lambdaonv}
\Lambda_{i,\lambda(h_i)}\Lambda_{i,-r}v = \Lambda_{i,\lambda(h_i)-r}v \quad\text{for all } i\in I, 0\le r\le \lambda(h_i).
\end{equation}
In other words, given $v, \lambda$, and $\omega_{i,r}$ as in the
proposition and setting
$$\gb\omega_{i}(u) = 1+\sum_{r=1}^{\lambda(h_i)} \omega_{i,r}u^r,$$
we want to show that
$$\Lambda_{i}^-(u)v = \gb\omega_i^-(u)v,$$
where for a polynomial $f(u) = \prod_j (1-a_ju)\in \mathbb F[u]$, we set $f^-(u) = \prod_j (1-a_j^{-1}u)$ (when convenient we shall also write $f=f^+$). The element $(\gb\omega_i)_{i\in I}$ is
called the Drinfeld polynomial of the $\ell$-highest-weight module
generated by $v$. We denote  by $\cal P^+_\mathbb F$ the multiplicative
monoid consisting of all $|I|$-tuples of the form
$\gb\omega = (\gb\omega_i)_{i\in I}$ where each $\gb\omega_i$ is a polynomial in
$\mathbb F[u]$ with constant term $1$.
The differential equations techniques used in \cite{cpweyl} for proving \eqref{e:Lambdaonv} do not work in positive characteristic. However, in light of Proposition \ref{ellhwrel}, it suffices to exhibit, for each $\gb\omega\in\cal P_\mathbb F^+$, one finite dimensional $\ell$-highest-weight module with $\ell$-highest weight $\gb\omega$ on which \eqref{e:Lambdaonv} is satisfied. This will be done in the next subsection.

We end this subsection introducing additional notation.  The  multiplicative
group corresponding to $\cal P_\mathbb F^+$ will be denoted by $\cal P_\mathbb F$. We let $\wt:\cal P_\mathbb F\to P$ be
the unique group homomorphism such that $\wt(\gb\omega) = \sum_{i\in
I} {\rm deg}(\gb\omega_i) \omega_i$ for all $\gb\omega\in\cal
P_\mathbb F^+$.
We also have an injective group homomorphism $\cal P_\mathbb F\to (U(\tlie h)_\mathbb F)^*$ given as follows. Any element $\gb\varpi\in\cal P_\mathbb F$ can be written uniquely as $\gb\omega\gb\pi^{-1}$, where $\gb\omega,\gb\pi\in\cal P_\mathbb F^+$ are such that $\gb\omega_i,\gb\pi_i$ are coprime for all $i\in I$. Its image  $\overline{\gb\varpi}\in(U(\tlie h)_\mathbb F)^*$ is defined  by
$$\overline{\gb\varpi}\left(\binom{h_i}{p^k}\right) = \binom{\wt(\gb\varpi)(h_i)}{p^k}, \qquad \overline{\gb\varpi}(\Lambda_i^\pm(u)) = \gb\varpi_i^\pm(u),$$
for all $k\in\mathbb Z_+$, and where $\gb\varpi_i^+=\gb\varpi_i$,
$\gb\varpi_i^-=\gb\omega_i^-(\gb\pi_i^-)^{-1}$. The second equality
is that of power series in $u$, obtained by expanding
$(\gb\pi_i^\pm)^{-1}$ as a product of geometric power series.  We
shall identify $\cal P_\mathbb F$ with its image in $(U(\tlie h)_\mathbb F)^*$ and
refer  to its elements as the integral $\ell$-weights. Similarly,
the elements in $\cal P_\mathbb F^+$ will be referred to as the dominant
integral $\ell$-weights.

\subsection{Classification of Irreducible Modules}\label{s:cirhla}
If $V$ is a finite-dimensional irreducible $U(\tlie g)_\mathbb F$-module,
proceeding as in the proof of Theorem \ref{cig}, we see that $V$ is
generated by a vector $v$ satisfying
$$(x_{\alpha,r}^+)^{(p^k)}v = 0,\  \ \binom{h_i}{p^k}v = \binom{\lambda(h_i)}{p^k}v,\ \ \Lambda_{i,r}v = \omega_{i,r}v,$$
for all $\alpha\in R^+, i\in I, r,k\in\mathbb Z, k\ge 0$ and some $\lambda\in P^+$, $\omega_{i,r}\in \mathbb F$. In particular, we have the following immediate corollary of Proposition \ref{ellhwrel}:

\begin{cor}\label{ilgehw}
Every finite-dimensional irreducible $U(\tlie g)_\mathbb F$-module is an
$\ell$-highest-weight module whose $\ell$-highest weight lies in
$\cal P_\mathbb F^+$.\hfill\qedsymbol
\end{cor}

We now introduce an important class of $U(\tlie g)_\mathbb F$-modules
called evaluation representations.

\begin{prop}\label{p:evmap}
For $a\in \mathbb F^{\times}$, there exists a surjective algebra homomorphism ${\rm hev}_a: U(\tlie g)_\mathbb F\to U(\lie g)_\mathbb F$ mapping $(x_{\alpha,r}^\pm)^{(k)}$ to $a^{rk}(x_{\alpha}^\pm)^{(k)}$. In particular, ${\rm hev}_a(\Lambda_{\alpha,r}) =(-a)^r\binom{h_\alpha}{|r|}$.
\end{prop}

We call ${\rm hev}_a$ the hyper evaluation map at $a$.

\begin{proof}
First observe that the formal evaluation map ${\rm ev}$ on $U(\tlie g)$ (see Lemma \ref{formalev}) sends $U(\tlie g)_\mathbb Z$ to $U(\lie g)_\mathbb Z\otimes \mathbb Z[t,t^{-1}]$. Hence, by reducing ${\rm ev}$ modulo $p$, we obtain the formal hyper evaluation map ${\rm hev}:U(\tlie g)_\mathbb F\to U(\lie g)_\mathbb F\otimes \mathbb F[t,t^{-1}]$. The statements of the proposition are now obvious (cf. definition of ${\rm ev}_a$ and \eqref{evLambda}).
\end{proof}

 Given any $U(\lie g)_\mathbb F$-module $V$,
let $V(a)$ be the pull-back of $V$ by ${\rm hev}_a$.  $V(a)$ is
called the evaluation representation with spectral parameter $a$
corresponding to $V$. For $a\in\mathbb F^\times$ and $\mu\in P$, let $\gb\omega_{\mu,a}$ be the element in $\cal P_\mathbb F$ whose $i$-th entry is $(1-au)^{\mu(h_i)}$, $i\in I$.
 If $V$ is a $U(\lie g)_\mathbb F$-highest-weight
module of highest weight $\lambda\in P^+$, it is easy to see that
$V(a)$ is an $\ell$-highest-weight module with Drinfeld
polynomial $\gb\omega_{\lambda,a}\in\cal P_\mathbb F^+$ and that the action of $\Lambda_i^-(u)$ on the $\ell$-highest vector is given by \eqref{e:Lambdaonv}. We shall denote
the evaluation representation by $V(\lambda,a)$ when $V=V(\lambda)$
and by $W(\lambda,a)$ when $V=W(\lambda)$.

 If $\lambda\in P_p^+$, it is easy to see that $V(p^k\lambda,a)$ is isomorphic to $V(\lambda,a^{p^k})^{\tilde\phi^{k}}$, where $\tilde\phi$ is the Frobenius homomorphism defined in section \ref{frsec}. Moreover, for any $\lambda\in P^+$, Theorem \ref{steinberg} implies
\begin{equation}\label{steinbergl}
V(\lambda,a)\cong \otimes_k V(p^k\lambda_k,a), \quad\text{where}\quad \lambda_k\in P^+_p\text{ are such that } \lambda = \sum_kp^k\lambda_k.
\end{equation}

We now prove the following version of Steinberg's Tensor Product Theorem for hyper loop algebras.

\begin{thm}\label{tpd} If $\mu_j\in P^+_p-\{0\}, a_j\in \mathbb F^{\times}$, and $l_j\in\mathbb Z_+, j=0,\cdots,n$, then $V=\otimes_j V(p^{l_j}\mu_j,a_j)$ is irreducible if and only if  $a_j\ne a_{j'}$ whenever $l_j=l_{j'}$.
\end{thm}

\begin{proof} The proof is a combination of the arguments used in Theorem \ref{steinberg} and \cite[Theorem 1.7]{cpnew}.
First consider the case $V=V(p^l\lambda,a)\otimes V(p^l\mu,b)$, where $\lambda,\mu\in P_p^+$, and let $v=\sum_j v_j\otimes w_j\in V$ be such that $w_j$ are linearly independent. Using \eqref{comnil} we get
$$(x_{\alpha,r}^+)^{(k)}v=\sum_j\sum_{l+m=k} a^{rl}b^{rm}\left((x_{\alpha}^+)^{(l)}v_j\right)\otimes \left((x_{\alpha}^+)^{(m)}w_j\right).$$
Hence, if $a=b$, this implies $(x_{\alpha,r}^+)^{(k)}v=a^{rk}
(x_{\alpha}^+)^{(k)}v$, and it follows that, if $v$ generates a
$U(\lie g)_\mathbb F$-submodule of $V$, it also generates a $U(\tlie
g)_\mathbb F$-submodule of $V$. This proves the ``only if'' part.

Conversely, for each $l\in\mathbb Z_+$, let $J_l = \{j: l_j=l\}$ and
$V_l = \otimes_{j\in J_l} V(p^l\mu_j,a_j)$, so that $V\cong
\otimes_{l:J_l\ne\emptyset} V_l$. Now observe that $V_l \cong
\left(\otimes_{j\in J_l}  V(\mu_j,a_j^{p^l})\right)^{\tilde\phi^l}$.
The same arguments used in \cite[Theorem 1.7]{cpnew} show that
$\otimes_{j\in J_l}  V(\mu_j,a_j^{p^l})$ is irreducible as $\tlie
g_\mathbb F$-module and, therefore, $V_l$ is irreducible as $U(\tlie
g)_\mathbb F$-module. Now let $\{l_1,\cdots, l_m\}=\{l:J_l\ne \emptyset\}$ and suppose
$l_1<l_2<\cdots<l_m$. Set $V'=\otimes_{j=2}^m V_{l_j}$ and let
$v=\sum_i w_i\otimes w'_i\in V_{l_1}\otimes V'$ such that $w_i'$ are
linearly independent. Then $(x_{\alpha,r}^+)^{(p^{l_1})}v =  \sum_i
\left((x_{\alpha,r}^+)^{(p^{l_1})}w_i\right)\otimes w'_i$ and we see
that $v$ is an $\ell$-highest-weight vector only if $v=v_1\otimes
v'$ where $v_1$ is an $\ell$-highest-weight vector for $V_{l_1}$,
since we already know that $V_{l_j}$ is irreducible. Proceeding
inductively, similarly to the proof of Theorem \ref{steinberg}, we
conclude that  $v$ must be a multiple of the tensor product of the
$\ell$-highest-weight vectors.
\end{proof}

As a corollary, we obtain the classification of the irreducible representations for $U(\tlie g)_\mathbb F$ (cf. \cite{cint,cpnew,cpweyl}). It is easy to see that every element $\gb\varpi\in \cal P_\mathbb F$ can be uniquely decomposed as $\gb\varpi = \prod_j \gb\omega_{\mu_j,a_j}$ for some $\mu_j\in P$ and $a_i\ne a_j$.

\begin{cor}\label{tpdirlg}\hfill
\begin{enumerate}
\item If $\gb\omega = \prod \gb\omega_{\lambda_j,a_j}\in \cal P_\mathbb F^+$ with $a_i\ne a_j, i\ne j$, and $\lambda_j=\sum_kp^k\lambda_{j,k}$ with $\lambda_{j,k}\in P_p^+$, then $V=\otimes_{j,k} V(p^k\lambda_{j,k},a_j)$ is an irreducible $U(\tlie g)_\mathbb F$-module with $\ell$-highest weight $\gb\omega$.
In particular, \eqref{e:Lambdaonv} holds for $V$.
\item  The isomorphism classes of irreducible
finite-dimensional $U(\tlie g)_\mathbb F$-modules are in one-to-one
correspondence with the elements of $\cal P_\mathbb F^+$.
\end{enumerate}
\end{cor}

\begin{proof}
It is immediate from Theorem \ref{tpd} that $V$ is irreducible and,
therefore, has an $\ell$-highest weight in $\cal P_\mathbb F^+$. In
order to see that this $\ell$-highest weight is $\gb\omega$, one
easily computes the action of $\Lambda_i^+(u)$ on the
$\ell$-highest-weight vector using \eqref{comcart} and observing
that each tensor factor is an evaluation representation. The proof
of \eqref{e:Lambdaonv} is completed in a similar way by computing
the action of $\Lambda_i^-(u)$ on the $\ell$-highest-weight vector.
This completes the proof of part (a) from which part (b) follows
immediately.
\end{proof}

Given $\gb\omega\in\cal P_\mathbb F^+$, let us denote by $V(\gb\omega)$ an irreducible $U(\tlie g)_\mathbb F$-module with $\ell$-highest weight $\gb\omega$. If $\gb\omega = \prod \gb\omega_{\lambda_j,a_j}\in \cal P_\mathbb F^+$ with $a_i\ne a_j, i\ne j$, it follows from \eqref{steinbergl} and the corollary above that $V(\gb\omega)\cong \otimes_j V(\lambda_{j},a_j)$.

Let us also record the following corollary.

\begin{cor}\label{ellw}
If $V$ is a finite-dimensional $U(\tlie g)_\mathbb F$-module, then
$V_{\gb\varpi}\ne 0$ only if $\gb\varpi\in\cal P_\mathbb F$ and $V_\mu=
\opl_{\gb\varpi:\wt(\gb\varpi)=\mu}^{} V_{\gb\varpi}$.
\end{cor}

\begin{proof}
It suffices to prove the claim for irreducible representations. Using \eqref{comcart} and the last corollary, it is sufficient to consider the irreducible evaluation representations $V=V(\lambda,a)$ with $\lambda\in P^+$. But in this case we have $V_\mu = V_{\gb\omega_{\mu,a}}$  (see also Proposition \ref{indepa} below and its corollary).
\end{proof}

We end the subsection computing the dual representation of a given irreducible one.
Let $V$ be a finite-dimensional $U(\tlie g)_\mathbb F$-module. Exactly as in the case of $U(\lie g)_\mathbb F$, we see that the dual vector space $V^*$ can be equipped with a $U(\tlie g)_\mathbb F$-module structure and $V^{**}$ is naturally isomorphic to $V$. Moreover, if $W$ is another finite-dimensional $U(\tlie g)_\mathbb F$-module, usual Hopf algebra arguments prove that we have a natural isomorphism of $U(\tlie g)_\mathbb F$-modules
\begin{equation}\label{dualten}
(V\otimes W)^* \cong W^*\otimes V^*.
\end{equation}

Given $\gb\varpi=\prod_j \gb\omega_{\mu_j,a_j}\in\cal P_\mathbb F, a_i\ne a_j$, set $\gb\varpi^* = \prod_j \gb\omega_{-w_0\mu_j,a_j}$.
We have:

\begin{prop}\label{p:dualaf}
Let $\gb\omega\in\cal P^+_\mathbb F$ and $V=V(\gb\omega)$. Then $V^*\cong V(\gb\omega^*)$.
\end{prop}

\begin{proof}
Due to Theorem \ref{tpd} and \eqref{dualten}, it suffices to consider the case $\gb\omega = \gb\omega_{\lambda,a}$ for some $a\in\mathbb F^\times$ and $\lambda\in P^+$. Since in this case $V$ is an evaluation representation,  every weight vector of $V$ is also an $\ell$-weight vector and  $V_\mu = V_{\gb\omega_{\mu,a}}$. Choose a basis for $V$ consisting of weight vectors. Then it is easy to see using \eqref{antipode} and \eqref{evLambda} that if $v$ is a basis element of weight $\mu$, then its dual vector $v^*$ is an $\ell$-weight vector of $\ell$-weight $\gb\omega_{w_0\mu,a}$. In particular, since $V^*\cong V(-w_0\lambda)$ as $U(\lie g)_\mathbb F$, we conclude that $V^*$ is the evaluation representation $V(-w_0\lambda,a)$.
\end{proof}

\subsection{The Weyl Modules}\label{wm}
We now study the universal finite-dimensional $\ell$-highest-weight $U(\tlie g)_\mathbb F$-modules motivated by \cite{cpweyl}.

\begin{defn}
Given $\gb\omega = (\gb\omega_i)_{i\in I}\in \cal P_\mathbb F^+$, let
$W(\gb\omega)$ be the $U(\tlie g)_\mathbb F$-module generated by a vector
$v$ satisfying
\begin{gather}
(x_{\alpha,r}^+)^{(p^k)}v = 0, \ \ \binom{h_i}{p^k}v = \binom{\wt(\gb\omega)(h_i)}{p^k}v, \ \ \Lambda_{i,\pm s} v = (\gb\omega^\pm_i(u))_s v,\\\notag \\ \label{fdimrel}
(x_{\alpha,r}^-)^{(l)}v = 0,
\end{gather}
for all $\alpha\in R^+, i\in I, k,l,r,s\in\mathbb Z, s,k\ge 0,
l>\wt(\gb\omega)(h_\alpha)$. Here, as before,
$(\gb\omega^\pm_i(u))_s$ means the coefficient of $u^s$. We call
$W(\gb\omega)$ the Weyl module with $\ell$-highest weight
$\gb\omega$.
\end{defn}

It follows from \eqref{comutxh} that
$$W(\gb\omega) = \opl_{\mu\le \wt(\gb\omega)}^{} W(\gb\omega)_\mu.$$
Standard arguments show:

\begin{prop}
$W(\gb\omega)$ has a unique maximal submodule and, hence, a unique
irreducible quotient. \qedsymbol
\end{prop}

In particular, $V(\gb\omega)$ is the irreducible quotient of
$W(\gb\omega)$,  $\gb\omega\in\cal P^+_\mathbb F$. Moreover, it follows from Proposition \ref{ellhwrel} that every finite-dimensional $\ell$-highest weight module of $\ell$-highest weight $\gb\omega$ is
isomorphic to a quotient of $W(\gb\omega)$. Hence, in order to complete the proof of the universality of $W(\gb\omega)$, it remains to show that it is finite dimensional. We begin with:

\begin{prop}If $W(\gb\omega)_\mu\ne 0$, then $W(\gb\omega)_{w\mu}\ne 0$ for all $w\in \cal W$. In particular, $W(\gb\omega)_\mu\ne0$ only if $w_0\wt(\gb\omega)\le \mu\le \wt(\gb\omega)$.
\end{prop}

\begin{proof}
Using an argument identical to the one used in characteristic zero, it follows from \eqref{fdimrel} that every vector $w\in W(\gb\omega)$ lies inside a finite-dimensional $U(\lie g)_\mathbb F$-submodule of $W(\gb\omega)$. Now all the claims follow from the corresponding results for finite-dimensional $U(\lie g)_\mathbb F$-modules.
\end{proof}

We are ready to prove:

\begin{thm}\label{wmfd}
$W(\gb\omega)$ is finite-dimensional for all $\gb\omega\in\cal P_\mathbb F^+$.
\end{thm}

This was proved in \cite{cpweyl}  for characteristic zero and for quantum groups, the
later under the assumption that $\lie g$ is simply laced (for non simply laced it follows from \cite{bn}).

\begin{proof}
Set $\lambda = \wt(\gb\omega)$ and let $v$ be an $\ell$-highest-weight vector of $W(\gb\omega)$.
It suffices to prove that $W(\gb\omega)$ is spanned by the elements
$$(x_{\beta_1,s_1}^-)^{(k_1)}\cdots (x_{\beta_m,s_m}^-)^{(k_m)}v,$$
with $m,s_j,k_j\in\mathbb Z_+, \beta_j\in R^+$ such that $s_j<\lambda(h_{\beta_j})$ and $\sum_j k_j\beta_j\le \lambda-w_0\lambda$. The last condition is immediate from  the previous proposition. Moreover,  the elements $(x_{\beta_1,s_1}^-)^{(k_1)}\cdots (x_{\beta_m,s_m}^-)^{(k_m)}v$ with no restriction on $s_j$ clearly span $W(\gb\omega)$.

Let $\cal R = R^+\times \mathbb Z\times\mathbb Z_+$ and $\Xi$ be the set of functions $\xi:\mathbb N\to \cal R$ given by $j\mapsto \xi_j=(\beta_j,s_j,k_j)$, such that $k_j=0$ for all $j$ sufficiently large. Let also $\Xi'$ be the subset of $\Xi$ consisting of the elements $\xi$ such that $0\le s_j<\lambda(h_{\beta_j})$. Given $\xi\in\Xi$ we associate an element $v_\xi\in W(\gb\omega)$ as above in the obvious way, i.e., if $k_j=0$ for $j>m$, then $v_\xi = (x_{\beta_1,s_1}^-)^{(k_1)}\cdots (x_{\beta_m,s_m}^-)^{(k_m)}v$. Define the degree of $\xi$ to be $d(\xi) = \sum_j k_j$ and the maximal exponent of $\xi$ to be $e(\xi)= \max\{k_j\}$.
Clearly $e(\xi)\le d(\xi)$ and $d(\xi)\ne 0$ implies $e(\xi)\ne 0$. Since there is nothing to be proved when $d(\xi)=0$ we assume from now on that  $d(\xi)>0$. Thus, let $\Xi_{d,e}$ be the subset of $\Xi$ consisting of those $\xi$ satisfying $d(\xi)=d$ and $e(\xi)=e$, and set $\Xi_d=\bigcup\limits_{1\le e\le d} \Xi_{d,e}$.

We prove by induction on $d$ and sub-induction on $e$ that if
$\xi\in\Xi_{d,e}$ is such that there exists $j$ with either
$s_{j}<0$ or $s_j\ge \lambda(h_{\beta_{j}})$, then $v_\xi$ is in the
span of vectors associated to elements in $\Xi'$. More precisely,
given $0<e\le d\in\mathbb N$, we assume, by induction hypothesis,
that this statement is true for every $\xi$ which belongs either to
$\Xi_{d,e'}$ with $e'<e$ or to $\Xi_{d'}$ with $d'<d$.

Observe that \eqref{basicrel} implies
\begin{equation}\label{basicrelv}
\left((X_{\beta;r,+}^-(u))^{(k-l)}\Lambda_{\beta}^+(u)\right)_kv = 0 \qquad \forall\ \beta\in R^+, k,l,r\in\mathbb Z, k>\lambda(h_\beta), 1\le l\le k.
\end{equation}
We split the proof in 2 cases according to whether $e=d$ or $e<d$.

When $e=d$, it follows that $v_\xi= (x_{\beta,s}^-)^{(e)}v$ for some $\beta\in R^+$ and $s\in\mathbb Z$. Suppose first that $e=1$ and let $l=\lambda(h_{\beta})$ and $k=l+1$ in \eqref{basicrelv} to get
\begin{equation}\label{basicrelv1}
(x^-_{\beta,r+1}\Lambda_{\beta,l} +x^-_{\beta,r+2}\Lambda_{\beta,l-1}+\cdots+x^-_{\beta,r+l+1})v=0.
\end{equation}
We consider the cases $s\ge l$ and $s<0$  separately  and prove the statement by  a further induction on $s$ and $|s|$, respectively. If $s\ge l$ this is easily done by setting $r=s-l-1$ in \eqref{basicrelv1}. Similarly, after observing that $\Lambda_{\beta,l}v\ne 0$, the case $s<0$ is dealt with by setting $r=s-1$ in \eqref{basicrelv1}. If $e>1$ let $l=e\lambda(h_{\beta})$ and $k=l+e$ in \eqref{basicrelv} to obtain
\begin{equation}\label{basicrelv2}
\sum_{n=0}^{\lambda(h_{\beta})} (x_{\beta,r+1+n}^-)^{(e)}\Lambda_{\beta,l-en}v + \text{ other terms} = 0,
\end{equation}
where the other terms belong to the span of elements $v_{\xi'}$ with $\xi'\in \Xi_{e,e'}$ for $e'<e$.
As before we argue by induction on $s$ and $|s|$ by setting  $r=s-1-\lambda(h_\beta)$ and $r=s-1$ in \eqref{basicrelv2}, respectively.

For the case $e<d$ we can assume, by inductions hypothesis, that $0\le s_j<\lambda(h_{\beta_j})$ for $j>1$. An easy application of Lemma \ref{sbrw} completes the argument in this case.
\end{proof}

\section{Reduction Modulo $p$}\label{redp}

\subsection{Introductory Remarks and Notation}\label{s:irmp}

In this section we start the theory of reduction modulo $p$ for
$U(\tlie g)_\mathbb K$-modules, where $\mathbb K$ is a field of
characteristic zero. In the case of $U(\lie g)_\mathbb K$ it
sufficed to prove the existence of admissible lattices for the
irreducible modules because the underlying abelian category was
semisimple. The category $\widetilde{\cal C}_{\mathbb K}$ is not
semisimple, so, even if it is possible to obtain a nice lattice
theory for all irreducible modules, one could not guarantee that all
of the objects in $\widetilde{\cal C}_{\mathbb K}$ would contain
such a lattice. In fact, even for irreducible modules the story is
more subtle than the one in the $U(\lie g)_\mathbb K$-case since the
evaluation maps ${\rm hev}_a$ do not preserve $\mathbb Z$-lattices
unless $a=\pm 1$. Still, we will prove that all the
$\ell$-highest-weight modules whose coefficients of their Drinfeld
polynomials are ``good'' with respect to $p$ can indeed be reduced
modulo $p$. In particular, it will follow  that every irreducible
$U(\tlie g)_\mathbb F$-module can be constructed as a quotient of a
module obtained by a reduction modulo $p$ process.

We consider two kinds of lattice theories. The first one is a
natural generalization of the one reviewed in Theorem \ref{elat} for
$U(\lie g)$.
Namely, in subsection
\ref{zlat}, we consider modules which contain finitely generated free
$\mathbb Z$-submodules which are invariant under the action of $U(\tlie g)_{\mathbb Z}$.
However, the modules $V^0(\lambda,a)$ with $a\in\mathbb
Z, a\ne \pm 1$, are easily seen not to contain such a lattice. Then,
in subsection \ref{slat}, we consider lattices over rings other than
$\mathbb Z$, namely, over torsion free discrete valuation rings. We think these
lattices are more suitable for studying reduction modulo $p$ in the
present context.

Let us fix some general notation to be used below. If $\mathbb A$ is any
commutative ring with identity,  let $\cal P_\mathbb A,\cal P^+_\mathbb A$ be defined in
the obvious way (cf. definition of $\cal P_\mathbb F$). Define also $\cal
P_\mathbb A^{++}$ as the subset of $\cal P_\mathbb A^+$ consisting of the elements
$\gb\omega$ such that the coefficient of the leading term of
$\gb\omega_i$ belongs to $\mathbb A^\times$ for all $i\in I$. Recall that
$\mathbb A$ is a discrete valuation ring if it is a local principal ideal
domain which is not a field and that its residue field is the
quotient of $\mathbb A$ by its unique maximal ideal. If $\mathbb A$ is a discrete
valuation ring with residue field $\mathbb F$, $a\in \mathbb A$, and $\gb\omega\in
\cal P_{\mathbb A}^{+}$, we let $\bar a$ and $\overline{\gb\omega}$ be the
images of $a$ in $\mathbb F$ and of $\gb\omega$ in $\cal P_{\mathbb F}^+$,
respectively. As before, $\mathbb F$ denotes an algebraically closed field
of characteristic $p>0$. We shall also denote by
$\overline{\gb\omega}$ the image of $\gb\omega\in\cal P^+_{\mathbb
Z}$ in $\cal P_\mathbb F^+$. We fix a torsion free discrete valuation ring $\mathbb A$ with residue field $\mathbb F$ (for instance, the ring of Witt vectors with coefficients in $\mathbb F$ \cite[Section II.5]{ser}) and denote by $\mathbb F^0$ the algebraic closure of the fraction field of $\mathbb A$.
Given $\gb\omega\in\cal P^+_{\mathbb F^0}$, we
denote by $W^0(\gb\omega)$ the corresponding $U(\tlie g)_{\mathbb F^0}$-Weyl module  \cite{cpweyl} and by $V^{0}(\gb\omega)$ its irreducible quotient.

\subsection{$\mathbb Z$-Lattices}\label{zlat}

\begin{defn}
If $V$ is a finite-dimensional $\mathbb F^0$-vector space we say that a finitely generated free $\mathbb Z$--submodule $L$ of $V$ is an ample lattice for $V$ if $L$ spans $V$ over $\mathbb F^0$. If the rank of $L$ is equal to the dimension of $V$, then we say $L$ is a lattice for $V$. If $V$ is a $U(\tlie g)_{\mathbb F^0}$-module, we say that an (ample) lattice for $V$ is admissible if $L$ is invariant under the action of
$U(\tlie g)_{\mathbb Z}$.
\end{defn}

If $L$ is an ample admissible lattice for a $U(\tlie g)_{\mathbb F^0}$--module $V$, we set $L_\mathbb F=L\otimes_{\mathbb Z}\mathbb F$. Thus, $L_\mathbb F$ is a $U(\tlie g)_\mathbb F$-module and ${\rm rank}(L)=\dim_\mathbb F(L_{\mathbb F})\ge \dim_{\mathbb F^0}(V)$.
It is trivial to see that the modules $V^0(\lambda,a)$ with $a\ne \pm 1$ do not contain a
finitely generated  $\mathbb Z$-submodule invariant under the action of $U(\tlie g)_{\mathbb Z}$. In fact, if $v$ is the $\ell$-highest-weight vector, then $(\Lambda_{i,\pm\lambda(h_i)})^kv = (-a)^{\pm k\lambda(h_i)}v$ is not a finitely generated $\mathbb Z$-module in that case.

\begin{prop}\label{zlatt}
Let  $V$ be a finite--dimensional $\ell$-highest weight $U(\tlie g)_{\mathbb F^0}$-module with
$\ell$-highest-weight $\gb\omega\in \cal P_{\mathbb Z}^{++}$ and
$\ell$-highest-weight vector $v$. Then $L=U(\tlie g)_{\mathbb Z}v$
is an   ample admissible lattice for $V$ and $L_\mathbb F$ is isomorphic to a
quotient of $W(\overline{\gb\omega})$. Moreover, if
$V=W^0(\gb\omega)$, then $L$ is a lattice.
\end{prop}

\begin{proof}
It is easy to see from \eqref{Lambda_alpha}, Lemma \ref{ht^snot0},
and \eqref{e:Lambdaonv} that $U(\tlie h)_{\mathbb Z}v=\mathbb Zv$
and, therefore, $L=U(\tlie n^-)_{\mathbb Z}v$. Also, $L$ is quite
clearly a torsion free $\mathbb Z$--submodule of $V$ which is
invariant under the action of $U(\tlie g)_{\mathbb Z}$. The proof of
Theorem \ref{wmfd} together with the hypothesis $\gb\omega\in \cal
P_{\mathbb Z}^{++}$ shows that  $L$ is a finitely generated $\mathbb
Z$-module which spans $V$ over ${\mathbb F^0}$ (the hypothesis
$\gb\omega\in \cal P_{\mathbb Z}^{++}$ is used to replace the remark
$\Lambda_{\beta,l}v\ne 0$ by $\Lambda_{\beta,l}v=av$ with
$a\in\mathbb Z^\times$). This completes the proof that $L$ is an
ample admissible lattice. Since the image of $v$ in $L_\mathbb F$ is
clearly an $\ell$-highest-weight vector with $\ell$-highest weight
$\overline{\gb\omega}$, the second statement follows immediately.
The last statement is clear since $L\otimes_{\mathbb Z}{\mathbb
F^0}$ is an $\ell$-highest-weight $U(\tlie g)$-module of
$\ell$-highest weight $\gb\omega$ and of dimension at least that of
$W^0(\gb\omega)$, thus $L\otimes_{\mathbb Z}{\mathbb F^0}\cong
W^0(\gb\omega)$.
\end{proof}

Clearly the only irreducible $U(\tlie g)_\mathbb F$-modules
which can be obtained as a quotient of some $L_\mathbb F$ where $L$ is as in the proposition
 are precisely those whose Drinfeld polynomials $\gb\omega$
lie in $\cal P^+_{\mathbb F_p}$ and the coefficient of the leading
term of $\gb\omega_i$ is $\pm 1$ for all $i\in I$. However, all of
the $\ell$-highest-weight $U(\tlie g)_\mathbb F$-modules $V$ whose Drinfeld
polynomial is of the form $\gb\omega_{\lambda,1}$ can be obtained in
this way. In the next section we will see that, for each $a\in\mathbb F^\times$, there exists an automorphism $\psi_a$ of $U(\tlie g)_\mathbb F$ determined by the assignment $(x_{\alpha,r}^\pm)^{(k)}\mapsto a^{rk}(x_{\alpha,r}^\pm)^{(k)}$ for all $\alpha\in R^+, k,r\in\mathbb Z, k>0$.
One can then show that the pull-back of such $V$ by $\psi_a$ is
an $\ell$-highest-weight module with Drinfeld polynomial
$\gb\omega_{\lambda,a}$. Hence, up to twisting by $\psi_a$, we obtain
all of the evaluation modules $V(\lambda,a)$. The other irreducible
modules are then obtained using tensor products.

\subsection{Lattices Over Discrete Valuation Rings}\label{slat}

We begin by giving a motivation for considering lattices over
discrete valuation rings.  Let $\mathbb P=\mathbb Z_{(p)}$ be the
localization of $\mathbb Z$ at $\mathbb Z-p\mathbb Z$ and $U(\tlie
g)_{\mathbb P}=U(\tlie g)_\mathbb Z\otimes_\mathbb Z \mathbb P$.
Then $\mathbb P$ is a torsion free discrete valuation ring with
residue field $\mathbb F_p$ and $U(\tlie g)_\mathbb P\otimes_\mathbb
P \mathbb F\cong U(\tlie g)_\mathbb F$. Let $a\in \mathbb P^\times$,
$v$ an $\ell$-highest-weight vector of $V=V^0(\lambda,a)$, and
$L=U(\tlie g)_{\mathbb P}v$. It is easy to see from
\eqref{Lambda_alpha},  \eqref{evLambda}, and Lemma \ref{ht^snot0}
that $U(\tlie h)_\mathbb Fv=\mathbb Pv$ and, therefore, $L=U(\tlie
n^-)_{\mathbb P}v = \mathbb P\left(U(\tlie n^-)_{\mathbb Z}v\right)=
\mathbb PL'$ where $L'=U(\lie n^-)_{\mathbb Z}v$ and $PL'$ is its
$\mathbb P$-span. Since $L'$ is the $\mathbb Z$-span of a basis for
$V$ by Theorem \ref{elat}, it follows that $L$ is the $\mathbb
P$-span of the same basis. Thus, setting $L_\mathbb F =
L\otimes_{\mathbb P} \mathbb F$, we obtain a $U(\tlie g)_\mathbb
F$-module isomorphic to $W(\lambda,\bar a)$, where $\bar a$ is the
image of $a$ in $\mathbb F_p$. This way we are able to obtain all
evaluation representations of the form $V(\lambda,b), b\in\mathbb
F_p$, as quotients of the reduction modulo $p$ of the irreducible
$U(\tlie g)$-modules $V^0(\lambda, a)$, where $a$ is such that $\bar
a=b$.  In order to obtain $V(\lambda,b)$ for all $b\in \mathbb F$ we
will have to use in place of $\mathbb P$ the bigger discrete
valuation ring $\mathbb A$ fixed in section \ref{s:irmp}. Recall
that $\mathbb F^0$ denotes the algebraic closure of the fraction
field of $\mathbb A$.

\begin{defn}\label{plat}
If $V$ is a finite-dimensional $\mathbb F^0$-vector space, we say
that a finitely generated free $\mathbb A$-submodule $L$ of $V$ is
an ample $\mathbb A$-lattice for $V$ if $L$ spans $V$ over $\mathbb
F^0$. If the rank of $L$ is equal to the dimension of $V$, then we
say $L$ is an $\mathbb A$-lattice for $V$. If $V$ is a $U(\tlie
g)_{\mathbb F^0}$-module, we say that an (ample) lattice for $V$ is
admissible if $L$ is invariant under the action of $U(\tlie
g)_{\mathbb A} = U(\tlie g)_{\mathbb Z}\otimes_{\mathbb Z}\mathbb
A$.
\end{defn}

 If $L$ is an ample $\mathbb A$-lattice for a $U(\tlie g)_{\mathbb F^0}$-module $V$,  we set $L_\mathbb F=L\otimes_{\mathbb A}\mathbb F$.
Then  $U(\tlie g)_\mathbb F\cong U(\tlie g)_\mathbb A\otimes_\mathbb A \mathbb F$ and $L_\mathbb F$ is a $U(\tlie
g)_\mathbb F$-module.  The next lemma is immediate.

\begin{lem}\label{tell}
Let $V$ and $W$ be finite-dimensional $U(\tlie g)_{\mathbb F^0}$-modules, $L$ and $M$ (ample) admissible lattices for
$V$ and $W$, respectively. Then $L\otimes_\mathbb A M$ is an (ample)
admissible lattice for $V\otimes W$ and $(L\otimes_\mathbb A M)_\mathbb F\cong
L_\mathbb F\otimes M_\mathbb F$ as $U(\tlie g)_\mathbb F$-modules.\hfill\qedsymbol
\end{lem}

\begin{thm}\label{latticeshg} Let $V$ be a finite-dimensional $U(\tlie g)_{\mathbb F^0}$-$\ell$-highest-weight module with  Drinfeld polynomial $\gb\omega\in\cal P_\mathbb A^{++}$ and $\ell$-highest-weight vector $v$. If $L = U(\tlie g)_{\mathbb A}v$ we have:
\begin{enumerate}
\item $L$ is an  ample admissible $\mathbb A$-lattice for $V$ and  $L_\mathbb F$ is isomorphic to a quotient of $W(\overline{\gb\omega})$.
\item If $V=W^0(\gb\omega)$, then $L$ is a lattice.
\item If  $V=V^0(\gb\omega)$ and $\gb\omega = \prod_{j=1}^m \gb\omega_{\lambda_j,a_j}$ with $\lambda_j\in P^+, a_j\in \mathbb A^\times, a_i\ne a_j$ when $i\ne j$, then $L$ is a lattice.
\end{enumerate}
\end{thm}

\begin{proof}
The proof of parts (a) and (b) are analogous to that of Proposition \ref{zlatt} with $\mathbb A$ in place of $\mathbb Z$.

We now prove (c). When $V$ is an
evaluation representation, i.e. when $m=1$, we proceed similarly to the motivational
discussion at the beginning of this subsection replacing $\mathbb P$
with $\mathbb A$, $U(\tlie g)$ with $U(\tlie g)_{\mathbb F^0}$, and
regarding $U(\tlie g)_\mathbb Z$ as embedded in $U(\tlie g)_{\mathbb
F^0}$. In the general case we have $V =
V^0(\lambda_1,a_1)\otimes \cdots\otimes V^0(\lambda_m,a_m)$. Let
$v_j$ be $\ell$-highest-weight vectors of $V^0(\lambda_j,a_j)$ so
that $v=v_1\otimes\cdots\otimes v_m$ and set $L' =
L_1\otimes_{\mathbb A} \cdots \otimes_{\mathbb A} L_m$, where $L_j =
U(\tlie g)_{\mathbb A}v_j$. Then $L_j$ are admissible lattices for
$V^0(\lambda_j,a_j)$ by the $m=1$ case and, by Lemma \ref{tell}, $L'$ is an admissible
lattice for $V$. It is clear from \eqref{comnil} that $L$ is an
$\mathbb A$-submodule of $L'$. Moreover, by part (a),  $L$ is a finitely generated free $\mathbb A$-module which spans $V$ and, hence, $L=L'$ since $\mathbb A$ is a principal ideal domain.
\end{proof}

The next  corollary  states that we have accomplished the task of
constructing all of the irreducible $U(\tlie g)_\mathbb F$-modules directly
as quotients of some $U(\tlie g)_{\mathbb F^0}$-modules by a reduction modulo $p$
process.

\begin{cor}\label{allirr}
For every $\gb\varpi\in\cal P_\mathbb F^+$ there exists $\gb\omega\in\cal P_\mathbb A^{++}$ such that
$\overline{\gb\omega}=\gb\varpi$ and
$V(\gb\varpi)$ is isomorphic to a quotient of $L_\mathbb F$,
where $L = U(\tlie g)_\mathbb Av$ and $v$ is an $\ell$-highest-weight
vector for $W^0(\gb\omega)$.
\end{cor}

\begin{proof}
Write $\gb\varpi = \prod_j \gb\omega_{\lambda,b_j}, b_j\in \mathbb F^\times, b_i\ne b_j$ for $i\ne j$, and let  $a_j\in\mathbb A^\times$ be lifts of $b_j$ to $\mathbb A$. Then the corollary follows from Theorem \ref{latticeshg} with  $\gb\omega=\prod_j\gb\omega_{\lambda_j,a_j}$.
\end{proof}

Let $\gb\omega, v,$ and $L$ be as in Theorem \ref{latticeshg},
suppose $V=W^0(\gb\omega)$, and write $\gb\omega = \prod_{j=1}^m
\gb\omega_{\lambda,a_j}$ with $a_i\ne a_{j}, i\ne j$, so that
$W^0(\gb\omega)\cong \otimes_j W^0(\gb\omega_{\lambda_j,a_j})$ (see \cite{cpweyl}).
Choose $\ell$-highest-weight vectors $v_j$ of
$W(\gb\omega_{\lambda_j,a_j})$ such that $v=v_1\otimes\cdots\otimes
v_m$ and set $L_j = U(\tlie g)_\mathbb Av_j$, $L'=L_1\otimes_\mathbb A\cdots \otimes_\mathbb A L_m$. As before, it follows from
\eqref{comnil} that $L\subseteq L'$.

\begin{con}\label{cp=0}
In the notation above we have:
\begin{enumerate}
\item $W(\overline{\gb\omega})\cong L_\mathbb F$.
\item If $\bar a_i\ne \bar a_j$ for $i\ne j$, then $L=L'$.
\end{enumerate}
\end{con}

Part (a) is the analogous statement of the conjecture in \cite{cpweyl} mentioned in the introduction of the paper. Notice that Theorem \ref{latticeshg} implies that $\dim_\mathbb F(L_\mathbb F)=\dim_{\mathbb F^0}(W^0(\gb\omega))$. Hence, for proving (a), it
suffices to prove that $\dim_\mathbb F(W(\overline{\gb\omega}))\le \dim_{\mathbb C}(W^0(\gb\omega))$.

Now part (b) is rather unusual since for $\mathbb Z$-lattices the appropriate analogous statement is false (as a counter-example one can take $\lie g=\lie{sl}_2, p\ne 2$, and $\gb\omega = (1-u)(1+u)$).
 Below we give an  example showing
that equality  can indeed happen when working with $\mathbb A$-lattices. This
is actually the main point behind the choice of working with
discrete valuation rings: they have plenty of units.  We have the following corollary of the Conjecture:

\begin{cor}\label{tpdwm}
Let $\lambda_j\in P^+$, $b_j\in\mathbb F^\times, j=1,\dots k$, be such that $b_i\ne b_j$ for $i\ne j$, and $\overline{\gb\omega} = \prod_j \gb\omega_{\lambda_j,b_j}$. Then:
\begin{enumerate}
\item $W(\overline{\gb\omega})\cong \otimes W(\gb\omega_{\lambda_j,b_j})$.
\item If $M_j$ is a quotient of $W(\gb\omega_{\lambda_j,b_j})$, $M=\otimes_j M_j$ is a quotient of $W(\overline{\gb\omega})$.
\end{enumerate}
\end{cor}

\begin{proof}
Let $a_j\in\mathbb F^0$ be such that $\bar a_j=b_j$. From part (a) of the conjecture we have $W(\overline{\gb\omega})\cong L_\mathbb F$ and from part (b) $L_\mathbb F = L'_\mathbb F$. Now Lemma \ref{tell} implies $L'_\mathbb F \cong \otimes_j (L_j)_\mathbb F$. Thus, applying part (a) of the conjecture to $(L_j)_\mathbb F$ we conclude part (a) of the corollary.

Once we have part (a), the proof of part (b) is standard. Namely, let $V_j$ be the kernel of the projection $W(\gb\omega_{\lambda_j,b_j}) \to M_j$.  Proceeding recursively on $j=1,\cdots,k$, we obtain short exact sequences
\begin{gather*}
0\to \left(\oti_{i=1}^{j-1} M_i\right)\oti_{}^{} V_j\oti_{}^{} \left(\oti_{i=j+1}^k W(\gb\omega_{\lambda_i,b_i})\right) \to \left(\oti_{i=1}^{j-1} M_i\right)\oti_{}^{}\left(\oti_{i=j}^k W(\gb\omega_{\lambda_i,b_i})\right) \to\\
\to \left(\oti_{i=1}^{j} M_i\right)\oti_{}^{} \left(\oti_{i=j+1}^k W(\gb\omega_{\lambda_i,b_i})\right) \to 0.
\end{gather*}
\end{proof}

The characteristic zero counterpart of part (a) of the corollary was
proved in \cite[Section 3]{cpweyl}. So far we did not manage to
adapt or compliment those techniques. By transferring the problem to
the setting of $\mathbb A$-lattices, we expect that other
characteristic zero arguments, e.g., as in  \cite{flmdw}, will lead
to a proof of part (b) of the  Conjecture.

Let us record the following proposition which follows immediately from Theorem \ref{elat}(a).

\begin{prop}\label{p:lattws}
Let $V$ be a finite-dimensional $U(\tlie g)_\mathbb F$-module.
Every additive subgroup of $V$ which is invariant under the action of
$U(\tlie g)_\mathbb A$ is the direct sum of its intersection with the
weight-spaces of $V$.\hfill\qedsymbol
\end{prop}

We now give the example showing that part (b) of Conjecture \ref{cp=0} in the setting of discrete valuation rings may hold. Let $\lie g = \lie{sl}_2$. Since $I$ is a
singleton, we shall drop the index referring to the roots and write
$x_r^\pm,h_r,$ and $\Lambda_r$ instead of $x_{1,r}^\pm$, etc., and
shall also identify $P$ with $\mathbb Z$. We will verify part (b) of
Conjecture \ref{cp=0} for the Weyl module $V=W^0((1-au)^2(1-bu))$
where $a,b\in \mathbb A^\times$ for some discrete valuation ring $\mathbb A$ such
that $\bar a\ne \bar b$. In particular $W^0((1-au)^2(1-bu))\cong
W^0((1-au)^2)\otimes W^0(1-bu)$. Let $v_0,w_0$ be $\ell$-highest
weight vectors of $W^0((1-au)^2)$ and $W^0(1-bu)$, respectively.

$W^0(1-bu)$ is isomorphic to the evaluation representation
$V^0(1,b)$. It is then easy to see that $x_s^-w_0 = b^sx_0^-w_0$ for
all $s\in\mathbb Z$. Thus, letting $w_1=x_0^-w_0$, the set
$\{w_0,w_1\}$ is an $\mathbb A$-basis for $L_2 = U(\tlie g)_\mathbb Aw_0$.

Now consider $W^0((1-au)^2)$ and let $L_1 = U(\tlie g)_\mathbb Av_0$. Since $\wt((1-au)^2)=2$, letting $k>2$ in \eqref{basicrel} we get
\begin{equation}\label{basicrele}
\left((X_{\alpha;s,+}^-(u))^{(k-l)}\Lambda_{\alpha}^+(u)\right)_kv_0 = 0 \qquad \forall\ l,s\in\mathbb Z,1\le l\le k.
\end{equation}
Setting $k=3, l=2$ above, we get
$(x_{s+1}^-\Lambda_{2}+x_{s+2}^-\Lambda_1+x_{s+3}^-)v_0=0$. Since
$\Lambda_2v_0=a^2v_0$ and $\Lambda_1v_0 = -2av_0$, one easily proves
inductively  that
\begin{equation}\label{basicrele1}
x_s^-v_0 = sa^{s-1}x_1^-v_0 - (s-1)a^sx_0^-v_0, \quad\text{for all}\quad s\in\mathbb Z.
\end{equation}
Let $v_1 = x_0^-v_0$ and $v_3 = x_1^-v_0$. Thus we see that
$\{v_1,v_3\}$ is an $\mathbb A$-basis for the zero-weight space of
$W^0((1-au)^2)\cap L_1$. Now setting $k=3,l=1$ in \eqref{basicrele}, we get
$(x_{s+1}^-)^{(2)}\Lambda_1v_0+x_{s+1}^-x_{s+2}^-v_0 = 0$. Setting $s=-1$ we get
\begin{equation}\label{basicrele3}
x_1^-x_0^-v_0 = 2a(x_0^-)^{(2)}v_0
\end{equation}
and setting $s=0$ we get
\begin{equation}
2a(x_1^-)^{(2)}v_0 = x_1^-x_2^-v_0.
\end{equation}
Now using \eqref{basicrele1} and then \eqref{basicrele3} on the right hand side of the last equation gives
\begin{equation}\label{basicrele4}
(x_1^-)^{(2)}v_0= a^2(x_0^-)^{(2)}v_0.
\end{equation}
Finally, using \eqref{basicrele1}, \eqref{basicrele3}, and \eqref{basicrele4} we get
\begin{equation}\label{basicrele2}
x_r^-x_s^-v_0 = 2a^{r+s}(x^-_0)^{(2)}v_0 \ \ \ \forall\ r,s\in\mathbb Z.
\end{equation}
Hence, $v_2=(x^-_0)^{(2)}v_0$ completes an $\mathbb A$-basis for $L_1$,
i.e., $L_1$ is the $\mathbb A$-span of $\{v_0,v_1,v_2,v_3\}$.

Clearly the set $A=\{v_i\otimes w_j: i=0,1,2,3 \text{ and } j=0,1\}$
is an $\mathbb A$-basis for $L'=L_1\otimes_\mathbb A L_2$. Since $L=U(\tlie
g)_\mathbb A(v_0\otimes w_0)\subseteq L'$, we are left to show that
$A\subseteq L$. Using \eqref{basicrele1}, \eqref{basicrele2} and
$x_s^-w_0 = b^sx_0^-w_0$ we compute:
\begin{align*}
x^-_0(v_0\otimes w_0) &= v_1\otimes w_0 + v_0\otimes w_1,\\
x^-_1(v_0\otimes w_0) &= v_3\otimes w_0 + bv_0\otimes w_1,\\
x^-_2(v_0\otimes w_0) &= 2av_3\otimes w_0-a^2v_1\otimes w_0 + b^2v_0\otimes w_1.
\end{align*}
Recording the coordinates of these vectors with respect to the basis
$\{v_1\otimes w_0, v_3\otimes w_0, v_0\otimes w_1\}$ of $L'\cap V_1$
we get the following matrix
$$\begin{bmatrix} 1&0&-a^2\\ 0&1&2a\\ 1&b&b^2\end{bmatrix}$$
whose determinant is $(a-b)^2$. Since $\bar a\ne \bar b$ iff $a-b\in
\mathbb A^\times$, we see that the vectors $x^-_0(v_0\otimes w_0)$,
$x^-_1(v_0\otimes w_0)$, and $x^-_2(v_0\otimes w_0)$ also form an
$\mathbb A$-basis for $L'\cap V_1$. Now we compute
\begin{align*}
(x^-_0)^{(2)}(v_0\otimes w_0) &= v_2\otimes w_0 + v_1\otimes w_1,\\
x^-_1x_0^-(v_0\otimes w_0) &= 2av_2\otimes w_0 + bv_1\otimes w_1 + v_3\otimes w_1,\\
(x^-_1)^{(2)}(v_0\otimes w_0) &= a^2v_2\otimes w_0 + bv_3\otimes w_1,
\end{align*}
and, recording the coordinates of these vectors with respect to the
basis $\{v_2\otimes w_0, v_1\otimes w_1, v_3\otimes w_1\}$ of
$L'\cap V_{-1}$ we get the matrix
$$\begin{bmatrix} 1&2a&a^2\\ 1&b&0\\ 0&1&b\end{bmatrix}.$$
The determinant  is again $(a-b)^2$ and we are done with this
weight-space as before. Finally one easily sees that
$(x_0^-)^{(3)}(v_0\otimes w_0) = v_2\otimes w_1$, showing that
$A\subseteq L$ as claimed.

Let us end this subsection proving the existence of the automorphisms $\psi_a$ as we promised at the end of section \ref{zlatt}.

\begin{prop}\label{p:homot}
For every $a\in\mathbb F^\times$ there exists an algebra automorphism $\psi_a:U(\tlie g)_\mathbb F\to U(\tlie g)_\mathbb F$ sending $(x_{\alpha,r}^\pm)^{(k)}$ to $a^{rk}(x_{\alpha,r}^\pm)^{(k)}$.
\end{prop}

\begin{proof}
Let $b$ be a lift of $a$ to $\mathbb A$ and let $\psi^0_b:U(\tlie g)_{\mathbb F^0}\to U(\tlie g)_{\mathbb F^0}$ be the algebra automorphism extending $t\mapsto bt$. It is easy to see that $\psi^0_b$ maps $U(\tlie g)_\mathbb A$ onto itself. Now let $\psi_a$ be the reduction modulo $p$ of the restriction of $\psi^0_b$ to $U(\tlie g)_\mathbb A$.
\end{proof}

\begin{rem}
Notice that the same kind of argument can be used to give an alternate proof of Proposition \ref{p:evmap} without using the formal evaluation map, but using the $\mathbb A$-form $U(\tlie g)_\mathbb A$.
\end{rem}

\subsection{Block Decomposition}\label{block}

We now assume Conjecture \ref{cp=0} in
order to obtain the block decomposition of the category of
finite-dimensional $U(\tlie g)_\mathbb F$-modules.
We begin with the following proposition on  the Jordan-H\"older series of Weyl modules.

\begin{prop}\label{indepa}
The $\ell$-weights of $W(\gb\omega_{\lambda,a})$ are of the form $\gb\omega_{\mu,a}$ with $\mu\in P$ such that $\mu\le \lambda$.
\end{prop}

\begin{proof}
Let $b$ be a lift of $a$ to $\mathbb A$, consider
$W^0(\gb\omega_{\lambda,b})$, and let $L=U(\tlie g)_{\mathbb F^0}v$
for some choice of $\ell$-highest-weight vector $v$ of
$W^0(\gb\omega_{\lambda,b})$. It is well-known that the
$\ell$-weights of $W^0(\gb\omega_{\lambda,b})$ are of the form
$\gb\omega_{\mu,b}$ with $\mu\in P$ such that $\mu\le \lambda$ (cf.
\cite[Proposition 3.3]{cmsc} and \cite{cmqc} for instance). In
particular, the weight spaces of $W^0(\gb\omega_{\lambda,b})$
coincide with its $\ell$-weight spaces and, therefore, using
Proposition \ref{p:lattws}, we conclude that $L$ is equal to its
intersection with the $\ell$-weight spaces of
$W^0(\gb\omega_{\lambda,b})$. Since, by Conjecture \ref{cp=0}(a),
$W(\gb\omega_{\lambda,a})$ is isomorphic to $L_\mathbb F$, the claim
of the proposition is now easily deduced.
\end{proof}

For each $a\in \mathbb F^{\times}$ and $i\in I$, set $\gb\omega_{i,a} = \gb\omega_{\omega_i,a}$ (the $\ell$-fundamental weights) and $\alpha_{i,a}(u) = \gb\omega_{\alpha_i,a}$ (the $\ell$-simple roots).
Let $\cal Q_\mathbb F$ (resp. $\cal Q_\mathbb F^+$) be the subgroup (resp. submonoid) of $\cal P_\mathbb F$
generated by all $\alpha_{i,a}(u)$. We call $\cal Q_\mathbb F$ the $\ell$-root lattice. We have the following Corollary of the preceding proposition together with Corollary \ref{tpdwm}.

\begin{cor}\label{wmsc}
If $V$ is a finite-dimensional $\ell$-highest-weight $U(\tlie
g)_\mathbb F$-module with $\ell$-highest weight $\gb\omega$,
$V_{\gb\varpi}\ne 0$ only if $\gb\varpi \in \gb\omega (\cal
Q_\mathbb F^+)^{-1}$.
\end{cor}

\begin{proof}
Proposition \ref{indepa} implies the result holds for $W(\gb\omega_{\lambda,a})$. Then we are done using Corollary \ref{tpdwm}. In fact \eqref{comcart} implies that the $\ell$-weights of the tensor product are products of the $\ell$-weights of each tensor factor (cf. \cite[Lemma 4.4]{cmqc}).
\end{proof}

\begin{defn}
A spectral character is a function $\chi: \mathbb F^{\times}\to P/Q$ with
finite support. Equipping the space of all spectral characters
$\Xi_\mathbb F$ with the usual abelian group structure, one sees that the
assignment $\gb\omega_{i,a}\mapsto \chi_{i,a}$, where $\chi_{i,a}(b)
= \delta_{a,b}\omega_i$, determines a group homomorphism $\cal
P_\mathbb F\to \Xi_\mathbb F, \gb\varpi\mapsto \chi_{\gb\varpi}$, with kernel $\cal
Q_\mathbb F$. We say that a $U(\tlie g)_\mathbb F$-module $V$ has spectral character
$\chi$ if  $\chi_{\gb\varpi} = \chi$ whenever $V_{\gb\varpi}\ne 0$.
Let $\widetilde{\cal C}_\chi$ be the category of all finite-dimensional
$U(\tlie g)_\mathbb F$-modules with spectral character $\chi$.
\end{defn}

We will denote by $\chi_{\mu,a}$ the spectral character corresponding to $\gb\omega_{\mu,a}, \mu\in P, a \in\mathbb F^\times$. We use additive notation for the group operation of $\Xi_\mathbb F$.

\begin{prop}\label{tpbinb}\hfill
\begin{enumerate}
\item For all $\gb\omega\in \cal P_\mathbb F^+$, $W(\gb\omega)\in\tilde{\cal C}_{\chi_{\gb\omega}}$.
\item $\widetilde{\cal C}_{\chi_1}\otimes\widetilde{\cal C}_{\chi_2}\subseteq\widetilde{\cal C}_{\chi_1+\chi_2}$ for all $\chi_1,\chi_2\in\Xi_\mathbb F$.
\item If $V\in \widetilde{\cal C}_\chi$ then $V^*\in \widetilde{\cal C}_{-\chi}$.
\end{enumerate}
\end{prop}

\begin{proof}
Parts (a) and (b) are immediate from Corollary \ref{wmsc} and its proof. Part (c) follows from Proposition \ref{p:dualaf}.
\end{proof}

Let $\widetilde{\cal C}_\mathbb F$ be the category of all finite-dimensional
$U(\tlie g)_\mathbb F$-modules. In the rest of the section we prove that the
block decomposition of $\widetilde{\cal C}_{\mathbb F}$ is described just as in the
characteristic zero case \cite{cmsc} and quantum group case
\cite{cmqc,em}. Namely:

\begin{thm}\label{bdlg}
The categories $\widetilde{\cal C}_\chi, \chi\in \Xi_\mathbb F$, are the blocks of $\widetilde{\cal C}_{\mathbb F}$.
\end{thm}

Once we have the statements of Propositions \ref{tpbinb} and \ref{p:dualfin} available, exactly the same arguments used in \cite[section 5]{cmsc} show that every indecomposable object from $\widetilde{\cal C}_{\mathbb F}$ belongs to some $\widetilde{\cal C}_\chi$, proving that we have the decomposition
$$\widetilde{\cal C}_{\mathbb F} = \opl_{\chi\in \Xi_\mathbb F}^{} \widetilde{\cal C}_\chi.$$
It remains to see that $\widetilde{\cal C}_\chi$ are indecomposable abelian subcategories. To do this it suffices to show that for any two given irreducible $U(\tlie g)_\mathbb F$-modules  $V$ and $W$ having the same spectral character, there exists a finite sequence of indecomposable objects $M_1, \cdots, M_k$ such that $V$ is a simple constituent of $M_1$, $W$ is a simple constituent of $M_k$ and, for every $j$, $M_j$ has a common simple constituent with $M_{j+1}$ (cf. \cite[Section 1]{em}). Let us begin with the case when $V=V(\lambda,a)$ and $W=V(\mu,b)$ for some $\lambda,\mu\in P^+$ and $a,b\in\mathbb F^\times$. Quite clearly $\chi_{\lambda,a}=\chi_{\mu,b}$ iff $\lambda-\mu\in Q$ and, if $\lambda\notin Q$, also $a=b$.

\begin{prop}\label{basicext}
Let $a\in\mathbb F^\times$ and suppose $\lambda,\mu\in P^+$ are such that ${\rm Hom}_{\lie g_{\mathbb F_0}}(\lie g_{\mathbb F_0}\otimes V^0(\lambda),V^0(\mu))\ne 0$ and $\lambda>\mu$. Then there exists a quotient $M$ of $W(\gb\omega_{\lambda,a})$ having $V(\mu,a)$ as simple constituent.
\end{prop}

\begin{proof}
Let $b$ be a lift of $a$ to $\mathbb A$. By \cite[Proposition 3.4]{cmsc}, there exists a non-split short exact sequence of $\tlie g_{\mathbb F_0}$-modules:
$$0\to V^0(\mu,b)\to M^0 \to V^0(\lambda,b)\to 0$$
for some $\ell$-highest-weight module $M^0$. From Theorem \ref{latticeshg}, there exists an ample admissible lattice $L$ for $M^0$ such that $M=L_\mathbb F$ is a quotient of $W(\gb\omega_{\lambda,a})$. It remains to show that there exists an $\ell$-highest-weight vector $v'$ for $V^0(\mu,b)$ in $M^0$ such $v'\in L$ and its image in $L_\mathbb F$ is non-zero. Thus, let $v$ be an $\ell$-highest-weight vector for $M^0$. From the proof of Theorem \ref{wmfd}, using that $b\in \mathbb A^\times$ as in the proof of Proposition \ref{zlatt}, we see that there exists an $\mathbb A$-basis for $L$ formed of vectors which are $\mathbb A$-linear combinations of elements of the form $(x_{\alpha_1,r_1}^-)^{(k_1)}\cdots (x_{\alpha_m,r_m}^-)^{(k_m)}v$. Let $v_1,\cdots, v_n$ be an $\mathbb A$-basis for $L_\mu$. Any $\ell$-highest-weight vector for $V^0(\mu,b)$ is a solution $\sum_{j=1}^n c_jv_j$, for some $c_j\in\mathbb F^0$, of the linear system
$$(x_{\alpha,r}^+)^{(k)}\left(\sum_{j=1}^n c_jv_j\right)=0$$
 for all $\alpha\in R^+, r\in\mathbb Z, k\in\mathbb Z_+$.  Since $L$ is admissible and the $\ell$-weights of $M^0$ are in $\cal P_\mathbb A$ (Proposition \ref{indepa}), it follows that there exists a solution with the $c_j$ lying in the field of fractions of $\mathbb A$. Since $\mathbb A$ is a unique factorization domain, it follows that we can choose $c_j$ in $\mathbb A$ such that the non-zero $c_j$ are coprime. This completes the proof.
\end{proof}

This proposition and Corollary \ref{tpdwm}(b) imply:

\begin{cor}\label{c:basicext}
Let $a=a_0,\lambda,\mu, M$ be as in Proposition \ref{basicext} and, for $j=1,\dots,k$, let $\nu_j\in P^+$ and $a_j\in \mathbb F^\times$ be such that $a_i\ne a_l$ for all $i,l = 0,\dots,k, i\ne l$. Then the $U(\tlie g)_\mathbb F$-module $M\otimes \left(\otimes_j V(\nu_j,a_j)\right)$ is $\ell$-highest-weight and has $V(\mu,a)\otimes \left(\otimes_j V(\nu_j,a_j)\right)$ and $V(\lambda,a)\otimes \left(\otimes_j V(\nu_j,a_j)\right)$ as simple constituents.\hfill\qedsymbol
\end{cor}

Now let $a\in\mathbb F^\times$ and $\lambda,\mu\in P^+$ be such that $\lambda-\mu\in Q-\{0\}$. Then by \cite[Proposition 1.2]{cmsc}, there exists a finite sequence $\lambda=\nu_1, \nu_2, \cdots, \nu_k=\mu$ such that $\nu_j\ne\nu_{j+1}$ and ${\rm Hom}_{\lie g_{\mathbb F_0}}(\lie g_{\mathbb F_0}\otimes V^0(\nu_j),V^0(\nu_{j+1}))\ne 0$. Since ${\rm Hom}_{\lie g_{\mathbb F_0}}(\lie g_{\mathbb F_0}\otimes V^0(\nu_j),V^0(\nu_{j+1})) = {\rm Hom}_{\lie g_{\mathbb F_0}}(\lie g_{\mathbb F_0}\otimes V^0(\nu_{j+1}),V^0(\nu_{j}))$, we conclude that there exists a sequence of $U(\tlie g)_\mathbb F$-$\ell$-highest-weight modules $M_j, j=1,\dots, k-1$, having both $V(\nu_j,a)$ and $V(\nu_{j+1},a)$ as simple constituents. From here it is quite clear how to complete the proof of Theorem \ref{bdlg} using the last corollary (cf. \cite[Section 4]{cmsc}).

\begin{rem}
We give an informal reasoning to justify why it should be expected
that the block decomposition of $\widetilde{\cal C}_\mathbb F$ is described
similarly to that of $\widetilde{\cal C}_{\mathbb F_0}$, contrary to what
happens with the block decompositions of $\cal C_\mathbb F$ and $\cal
C_{\mathbb F_0}$ (the categories of finite-dimensional representations
for $U(\lie g)_\mathbb F$ and $U(\lie g)_{F_0}$, respectively). While the blocks
of $\cal C_{\mathbb F_0}$ are as small as possible ($\cal C_{\mathbb F_0}$ is a semisimple category), the blocks of $\tilde{\cal C}_{\mathbb F_0}$ are as large as one can expect (for instance, when
$P/Q$ is trivial, $\widetilde{\cal C}_{\mathbb F_0}$ is itself an
indecomposable abelian category). Hence, while the blocks of $\cal
C_\mathbb F$ have space to become ``larger'' (and they indeed become so, but
still not as large as possible \cite[Chapter II.7]{janb}), that is
not the case for $\widetilde{\cal C}_\mathbb F$.
\end{rem}

\bibliographystyle{amsplain}

\end{document}